\documentclass{amsart}
\usepackage{graphicx}
\usepackage[english]{babel}
\usepackage{amsmath,amssymb}
\usepackage{amsthm}
\usepackage{todonotes}
\usepackage{tikz}
\usetikzlibrary{decorations.pathreplacing,calligraphy}
\usepackage{hyperref}
\usepackage{enumitem}
\theoremstyle{definition}
\newtheorem{defin}{Definition}[section]
\newtheorem{notat}[defin]{Notation}

\theoremstyle{plain}
\newtheorem{teo}[defin]{Theorem}
\newtheorem{prop}[defin]{Proposition}
\newtheorem{lemma}[defin]{Lemma}

\newtheorem{cor}[defin]{Corollary}

\theoremstyle{definition}
\newtheorem{remark}[defin]{Remark}

\newcommand{\RR}{\mathbb{R}}
\newcommand{\NN}{\mathbb{N}}
\newcommand{\ZZ}{\mathbb{Z}}
\newcommand{\QQ}{\mathbb{Q}}
\newcommand{\dom}{\operatorname{dom}}
\newcommand{\ran}{\operatorname{ran}}
\newcommand{\baire}{\mathbb{N}^{\mathbb{N}}}
\newcommand{\diam}{\operatorname{diam}}
\newcommand{\FF}{\mathcal{F}}
\newcommand{\CUBE}{\mathrm{CUBE}}

\usepackage{color} %used for NOTES
\newcommand{\guido}[2] %Eric's blue comments % 1. Arg.: no Text, 2. Arg.: highlighted Text
   {\textcolor{red}{#2}}

\newcommand{\andrea}[2] %Eric's blue comments % 1. Arg.: no Text, 2. Arg.: highlighted Text
   {\textcolor{blue}{#2}}

\title{Computability of a Whitney Extension}
\author[Brun]{Andrea Brun}
\address{Dipartimento di Scienze Matematiche, Informatiche e Fisiche, Universit\`{a} di Udine, Italy}%
\email{andrea.brun.tv@gmail.com}
\author[Gherardi]{Guido Gherardi}
\address{Dipartimento delle Arti, Universit\`{a} di Bologna, Italy; ORCID iD 0000-0002-3382-2292}
  \email{guido.gherardi@unibo.it}
\author[Marcone]{Alberto Marcone}
\address{Dipartimento di Scienze Matematiche, Informatiche e Fisiche, Universit\`{a} di Udine, Italy; ORCID iD 0000-0001-8356-0086}
\email{alberto.marcone@uniud.it}
\date{\today}

\thanks{The first author started working on this project for his master thesis, discussed at University of Udine in March 2024. The second and third author were partially supported by the  Italian PRIN 2022 ``Models, sets and classifications'', prot.\ 2022TECZJA, funded by the European Union - Next Generation EU. Marcone is a member of INdAM-GNSAGA}

\subjclass{Primary 03D78; Secondary 26E10, 58C25, 54C20, 54C30}
\keywords{Whitney Extension Theorems, Computable Analysis; Type-$2$ Theory of Effectivity}
\begin{document}

\begin{abstract}
We prove the computability of a version of Whitney Extension, when the input is suitably represented.
More specifically, if $F \subseteq \RR^n$ is a closed set represented so that the distance function $x \mapsto d(x,F)$ can be computed, and $(f^{(\bar{k})})_{|\bar{k}| \le m}$ is a Whitney jet of order $m$ on $F$, then we can compute $g \in C^{m}(\RR^n)$ such that $g$ and its partial derivatives coincide on $F$ with the corresponding functions of $(f^{(\bar{k})})_{|\bar{k}| \le m}$.
\end{abstract}

\maketitle

\section{Introduction}

In this paper, we investigate the computability of a version of Whitney Extension.
A similar investigation was already carried out by Charles Fefferman in \cite{Fefferman:2009} using the real RAM approach.
We tackle a similar question in the framework of Type-2 Theory of Effectivity approach to computability.\smallskip

In the seminal paper \cite{whitney:1934} Hassler Whitney proved a result, the so-called Whitney Extension Theorem, that gave rise to a line of research which is still active today (see e.g.\ \cite{ONeill:2025}).
The goal, roughly speaking, is to obtain generalizations to the case of differentiable functions of the Urysohn-Tietze Extension Theorem.
The version we deal with states that a partial continuous function $f$, defined on a closed subset $F$ of $\RR^n$ and differentiable up to order $m$ for some $m \ge 0$, can be extended to a total continuous function $g$ with the same degree of differentiability and such that its partial derivatives coincide with those of $f$ on $F$.
There is a subtlety here that requires more explanation.
%Since we cannot talk about differentiation of functions defined on closed sets, we must first clarify what we mean when we say that a function is differentiable on a closed subset of $\RR^n$.
%To do this, we employ the concept of \emph{Whitney jet} on a closed set $F$, which is a finite set of continuous functions defined on $F$ that satisfies Taylor's Theorem, and hence behaves like partial derivatives of each other (further details in Section 	\ref{Section_WETm}).
It must be in fact clarified in what sense one can say that a function is differentiable on a closed subset of $\RR^n$ and in particular on its boundary points.
For this purpose, it is customary to use \emph{Whitney jets}: these are finite sets of continuous functions defined on the same closed $F$ and  satisfying Taylor's property on $F$, i.e., they behave like partial derivatives of the leading function $f$ in the jet, seen as a derivative of degree 0 of itself (further details in Section 	\ref{Section_WETm}).
Being a Whitney jet is obviously a necessary condition for the existence of a differentiable extension of the set of functions, and in particular of $f$.
The version of Whitney Extension that we consider states that this condition is also sufficient.
We can thus restate the theorem more precisely as ``Given a closed set $F$ and a Whitney jet of order $m$ on $F$, there exists a total continuous function $g \in C^{m}(\RR^n)$ such that both $g$ and its partial derivatives coincide on $F$ with the corresponding functions of the jet''.
Moreover, once $F$ is fixed, the map from the Whitney jet to $g$ is linear.\smallskip

This version of Whitney Extension provides a typical example of what in Computable Analysis is called a \emph{mathematical problem}.
In this context, a mathematical problem is a true statement of the form $\forall x \, (\Phi (x) \to \exists y \,\Psi(x,y))$. An \emph{instance} of the problem is any $x$ satisfying $\Phi(x)$, and a \emph{solution} for this instance is some $y$ such that $\Psi(x,y)$.
A mathematical problem can be naturally expressed by a multi-valued function, which associates an element $x$ from the domain $\{x \mid \Phi(x) \}$ to the non-empty set of solutions $\{ y \mid \Psi(x,y) \}$ for the instance $x$.
A mathematical problem is computable if the associated multi-valued function is computable in the sense of Computable Analysis,  that is, by means of Type-2 Turing machines.\smallskip

In the case we deal with, an instance is a pair made of a closed set $F \subseteq \RR^n$ and a Whitney jet of a given order $m\geq 0$ defined on $F$, while a solution for this instance is a total differentiable function $g : \RR^n \to \RR$ that extends the leading function $f$ and whose partial derivatives up to degree $m$ coincide on $F$ with those of $f$.
%We call $\mathrm{WET}$ the multi-valued function associated to the Whitney Extension Theorem.

To analyze the computable properties of this mathematical problem, we look at the classical proof given by Stein in \cite[Chapter VI]{stein:1970} and examine it in detail.

Along these lines, our proof constitutes of two steps: first of all we consider only the property of continuity and obtain a \emph{First Computable Whitney Extension Theorem}, which also provides a new version of Computable Urysohn-Tietze Theorem (with the extra-condition of  linearity).
Then the condition of differentiability is further considered, obtaining our full \emph{Computable Whitney Extension Theorem}.
%In this way it is also proved how the Whitney Extension Theorem can be really seen as an extension of the Urysohn-Tietze Theorem.
%, modifying it appropriately to show that the resulting multi-valued functions is computable.

Klaus Weihrauch had already obtained a Computable Urysohn-Tietze Theorem by showing in \cite{WeiTietze} that the classical proof contained in \cite{engelking:1989} can be turned into an algorithm, so it is “essentially” computable\footnote{For a comparison with the Reverse Mathematics approach, one can refer to \cite{shafer:2016}, where several version of the Tietze Extension Theorem have been proved to be equivalent to $\mathsf{WKL}_0$ over $\mathsf{RCA_0}$.}
%It is then natural to pursue the same kind of investigation for the Whitney Extension Theorem and its proofs.
His result can be seen as an instance of a more general research line, with also philosophical significance, evaluating to what extent classical proofs of statements giving rise to mathematical problems have an algorithmic nature.

Some such proofs present indeed an intuitive computational flavor, but, so to say, some computable steps are hidden under the surface. In such cases, the proofs of computability of the corresponding multi-valued functions can retrace those classical proofs and show in a rigorous way the computational content that was only sketched in them. Nevertheless, removing the rind to get to the computational pith might contain non-trivial steps, and it depends first of all on the correct choice of suitable translations of the classical concepts into computational notions.

Stein's proof of a Whitney Extension result has a clear computational \lq\lq flavor\rq\rq, but to what extent is it really \lq\lq computable\rq\rq?
That is, can it be turned into an algorithm?
And, if so, at what cost?
%Or, in other words, what is the price of making it computable?

The proof at hand is based on four main ingredients: 1) a family of cubes tiling the complement of the closed domain $F$; 2) a smooth partition of unity defined on $\RR^n \setminus F$; 3) the use of projection points on $F$ of the cubes in the tiling\footnote{This step was the motivation behind \cite{proj:2019}.}; 4) a formula, based on the first three ingredients, explicitly defining the total extension $g$ of the differentiable $f$ outside its domain $F$.

Although these ingredients of the proof might appear at first sight computable, a closer look shows that this is not the case. In particular,
i) Stein's cube tiling cannot be computed from $F$;
ii) the problem of establishing which are the cubes of a tiling containing a given point is not recursively solvable;
iii) a partition of the unit made of functions that are not only smooth, but also \emph{computably smooth} is needed;
iv) projection points cannot be found computably (\cite{proj:2019});
v) Stein's definition of $g$ is based on a case distinction depending on  whether a point $x$ is in $F$ (in this case $g(x)=f(x)$), or not (in this case $g(x)$ is determined by Stein's formula); however this distinction is not recursive.

We will see how these problems can be overcome, provided that enough information is given about the closed domain $F$.
Points iv) and v) above are particularly challenging.

Regarding iv), Stein himself indirectly suggests a possible solution, as he observes, with respect to a generic cube $Q$ of the tiling, that \lq\lq in fact any choice of a point in $F$ with the property that the distance of that point to $Q$ is comparable to the distance of $Q$ from $F$ would do as well\rq\rq\ \cite[p.\ 172]{stein:1970}.
Nevertheless, he states that considering the projection points is simpler.
Well, it is the notion of \lq\lq simpler\rq\rq\ that needs now some more care. Conceptually, it is indeed simpler to consider points defined in a simple way, such as \lq\lq points of minimum distance\rq\rq, no matter how much this amounts to.
But, computably, what can be computed from the input is simpler than what is not. By this reason, we are moved rather to selects points that lie at a distance that is, as Stein says, \lq\lq comparable to the distance of $Q$ from $F$\rq\rq.
Points of this kind can be in fact computably found, provided that we have sufficient information on the closed $F$, once a quantitative notion of \emph{comparability} has been suitably determined.

The treatment of point v) will instead require the determination of a strategy based on the choice of suitable real constants assuring that the output computed at some stage $i$
will approximates $g(x)$ up to an error of $2^{-i}$, no matter whether $x$ is in $F$ or not, a question that in general we might not be able to answer at that stage.

It turns out that to solve the problems outlined above we can use the so-called \emph{total} information about $F$ (Definition \ref{Def_ClosedSetRepresentations}); moreover even the determination of the real constants involved in the solution of v) might appear sometimes tricky: however a selection can be made so that Stein's proof can be turned into an algorithm.

The main result of the paper is thus Theorem \ref{Theorem_WETm}, stating that the version of the Whitney Extension problem we consider is indeed computable when the original domain $F$ is given with total information. This fulfills a promise made in \S6 of \cite{proj:2019}, where the ideas of the present paper were briefly sketched.\smallskip

We now outline the organization of the paper.
Section \ref{Section_ComputableAnalysis} recalls the basic notions of Computable Analysis and sets the stage for analyzing Whitney Extension Theorem.
In Section \ref{Section_DecompositionCubes} we show how to obtain a computable tiling of the complement of a closed $F \subseteq \RR^n$ given with total information.
Section \ref{Section_PartitionUnit} deals with the partition of unity based on the tiling from the previous section.
In Section \ref{Section_FirstWET} we prove our First Computable Whitney Extension Theorem, while Section \ref{Section_WETm} establishes our full Computable Whitney Extension Theorem.
The final Appendix \ref{Appendix_A} includes the detailed proof of Proposition \ref{Prop_Derivativesvarphi}, a technical result providing bounds for the derivatives of the functions involved in the partition of unity.

\section{Computable Analysis}\label{Section_ComputableAnalysis}

In this section we review some fundamental notions of Computable Analysis as outlined in the Type-2 Theory of Effectivity (TTE).
We refer the reader to \cite{weihrauch:2000} for a more detailed introduction to the subject.
We additionally introduce some tools that are needed for the analysis of Whitney Extension Theorem.

%Ordinary computability theory first introduces computable partial functions $f : \subseteq \NN \to \NN$ explicitly, for example, by means of Turing machines. For defining computability on other sets $S$ (rational numbers, finite graphs, etc.) numbers are used as ``names'' or ``codes'' of elements of $S$.
%While a machine still transfers natural numbers to natural numbers, the user interprets these numbers as names of elements from the set $S$. Since the set $\NN$ is countable, it is not sufficient as set of names for uncountable sets such as $\RR$.

The idea underlying TTE is that there are types of mathematical objects, such as real numbers, that require an infinite amount of information to be determined in a precise way.
It is therefore natural to employ machines designed to compute with infinite sequences of natural numbers (rather than finite), which are called \emph{Type-2 Machines}:

%\subsection{Type-2 Machines}

%We assume that Turing machines have only one input tape and one work tape.
%This is equivalent to using Turing machines with multiple inputs and work tapes (see \cite{papadimitriou:1993} for further details).

\begin{defin}[Type-$2$ Turing machine]
    A \emph{Type-$2$ Turing machine} is a Turing machine with three tapes: an input tape, a work tape and an output tape.
    %On the input tape, the only operations allowed are moving to the right and reading, on the work tape all operations are allowed, and on the output tape the only operation allowed are writing and move to the right.
    %In other words, the input tape is one-way and read-only, while the output tape is one-way and write-only.
    The input tape is one-way and read-only, the output tape is one-way and write-only, whereas all the ordinary Turing machine operations (reading, writing, erasing) are allowed on the two-way work tape.
\end{defin}

%The restrictions imposed to Type-2 machines are designed to deny the possibility of making corrections on the output tape.	

We say then that a partial function $f:\subseteq \baire\to\baire$ is \emph{computable} if there is a Type-2 Turing machine writing $f(p)\in\baire$ on the output tape for any $p\in \dom(f)$ given on the input tape (we use  \lq\lq$\subseteq$\rq\rq for partial functions to mean that the domain is a subset of the source space).

%The restriction to one-way output guarantees that any initial segment of $f(p)$ computed at a certain stage cannot be erased in the future and, therefore, is final.
The restrictions concerning the output tape guarantee that any initial segment of $f(p)$ computed at a certain stage cannot be erased in the future.
In contrast, models with two-way output tapes and allowing erasures do not guarantee the correctness of partial results, and therefore are not  useful in practice.

The space $\baire$ of infinite sequences over $\NN$ is endowed with the usual \emph{Baire topology}, that is, the topology generated by the basic open sets
\begin{align*}
    w \baire := \left\{ p \in \baire \mid w \text{ is an initial prefix of }p \right\},
\end{align*}
where $w$ belongs to the set $\NN^*$ of all finite sequences over $\NN$.
The Baire space is homeomorphic to its own finite and countable products via the usual Cantor pairing functions and projections.

%The topological space $(\baire, \tau_B)$ is called \emph{Baire space}.
% open set in the Baire space has the form
%\begin{align*}
    %W \baire := \bigcup_{u \in W} u \baire,
%\end{align*}
%where $W$ is an arbitrary subset of $\NN^*$. The Baire space is second-countable, separable and zero-dimensional, since the basic open sets $u \baire$ are clopen.

%The Baire topology is induced by the metric $d_{B} : \baire \times \baire \to \RR_{\ge 0}$ defined by
%\begin{align*}
    %d_B(p,q) :=
    %\begin{cases}
        %0 & \text{ if }p = q, \\
        %2^{-\min\{i \in \NN \, \mid \, p(i) \ne q(i) \}} & \text{ otherwise.}
    %\end{cases}
%\end{align*}
%The metric $d_B$ is complete, meaning that every Cauchy sequence has a limit.

It is well known that all partial computable funtions $f:\subseteq\baire\to  \baire$ are continuous with respect to the relative topology induced by the Baire topology on $\dom(f)$, and that partial computable functions are closed under composition.

\subsection{Represented Spaces}\label{Section_RepresentedSpaces}
To extend the notion of computability to functions between spaces different from $\baire$ we need the notion of \emph{representation}.
%A representation equips the objects of a given space with names, giving rise to the concept of represented space. Computable functions between represented spaces are those which are realized by a computable function on the names.
A representation provides a naming system to encode objects of a given space (set) by using infinite sequences of natural numbers.

\begin{defin}[Represented space]
    A \emph{represented space} $(X, \delta)$ is a set $X$ together with a surjective map $\delta : \subseteq \baire \to X$, which is said to be a \emph{representation} of $X$.
    If $\delta (p) = x$, we call $p$ a \emph{$\delta$-name} for $x$ (or simply a \emph{name}, when the representation is clear from the context).
\end{defin}

%Notice that in the definition above, there is no requirement for $\delta$ to be injective.

An element $x$ of a represented space $(X, \delta)$ is called \emph{computable} if it has a computable name $p \in \baire$.
%, that is, a name $p$ that is output by some Type-2 machine on input $0^\mathbb N$.

Given represented spaces $X_1,X_2,...$ with corresponding representations $\delta_1,\delta_2,...$, the representations $\delta_1\times...\times\delta_n$ and $\prod_{i\in \mathbb N}\delta_i$ of the spaces $X_1\times...\times X_n$ and $\prod_{i\in \mathbb N}X_i$, respectively, can be obtained in an obvious way by using Cantor pairing functions $\langle \cdot, \dots, \cdot \rangle$ and $\langle \cdot, \dots \rangle$ on $\mathbb N^\mathbb N$.
Moreover, we can use pairing functions to combine different representations for a single space $X$:

%Similarly, a \emph{computable sequence} $(x_i)_i$ in $(X, \delta)$ is a sequence of points such that there exists a computable $p \in \baire$ with $\delta ( \pi_{\infty,i}(p) ) = x_i$ for all $i \in \NN$.

%\subsubsection{Construction of new representations}

%Now we present canonical methods of constructing new representations from given ones.

%\begin{defin}
    %Let $(X_1, \delta_1), \dots, (X_k, \delta_k)$ be represented spaces.
    %We define the representation $[\delta_1, \dots, \delta_k]$ of $X_1 \times \dots \times X_k$ by
    %\begin{align*}
        %[\delta_1, \dots, \delta_k] \langle p_1, \dots, p_k \rangle := (\delta_1 (p_1), \dots, \delta_k (p_k)).
    %\end{align*}
    %Let now $(X_i, \delta_i)$ be a represented space for every $i \in \NN$. We define the representation $[\delta_0, \dots]$ of $\prod_{i \in \NN} X_i$ by
    %\begin{align*}
     %   [\delta_0, \delta_1, \dots] \langle p_0, p_1, \dots \rangle := (\delta_0 (p_0), \delta_1 (p_1), \dots).
   % \end{align*}
    %We use the abbreviations $[\delta]^k$ for $[\delta, \dots, \delta]$, and $[\delta]^\omega$ for $[\delta, \dots]$.
%\end{defin}

%Therefore $p = \langle p_1, p_2 \rangle$ is a $[\delta_1, \delta_2]$-name for the pair $(x_1, x_2)$ iff $p_1$ is a $\delta_1$-name for $x_1$ and $p_2$ is a $\delta_2$-name for $x_2$.

\begin{defin}\label{Def_RepresentationsWedge}
For all $i \in \{1, \dots, k\}$ let $\delta_i : \subseteq \baire \to X$ be a representation of the set $X$.
We define the representation $\delta_1 \wedge \dots \wedge \delta_k$ of the set $X$ by
        \begin{align*}
            \delta_1 \wedge \dots \wedge \delta_k \langle p_1, \dots, p_k \rangle = x \iff \delta_1(p_1) = \dots = \delta_k (p_k) = x.
        \end{align*}
\end{defin}

\subsubsection{Computability of multi-valued functions between represented spaces}

Multi-valued functions $f : \subseteq X \rightrightarrows Y$, also called \emph{problems}, are just relations $f \subseteq X \times Y$, with a suitable notion of composition.
Let $\dom(f)$ the set of all \emph{admissible instances} of the problem $f$, that is $\dom(f) = \{ x \in X \mid \exists y \: (x,y) \in f \}$.
Given $x \in \dom(f)$, we define $f(x)$ as the set of all possible \emph{results} of the problem for $x$, i.e.\ $f(x) := \{ y \in Y \mid (x,y) \in f \}$.
%We call \emph{mathematical problem} (or simply \emph{problem}) any given partial multi-valued function $f : \subseteq X \rightrightarrows Y$ between represented spaces $X$ and $Y$.
The \emph{composition} of two problems $f:\subseteq X\rightrightarrows Y$ and $g:\subseteq Y\rightrightarrows Z$ is the problem $g\circ f:\subseteq X\rightrightarrows Z$ with $\dom(g\circ f):=\{x\in \dom(f) \mid f(x)\subseteq \dom(g)\}$ and $g\circ f(x):=\{z\in Z\mid \exists y\in f(x)\, z\in g(y)\}$.

%The set $X$ is called \emph{source} $f$ and $Y$ is called \emph{target} of $f$.

%		\begin{example}[\cite{brattka:2021}, Zero problem]
%			Let $X$ be a topological space and let $C(X)$ be the set of continuous functions $f : X \to \RR$. The \emph{zero problem} denoted as $\mathrm{Z}_X : \subseteq C(X) \rightrightarrows X$ is defined by the mapping $f \mapsto f^{-1}(0)$. In simpler terms, the zero problem is the problem of finding solutions $x \in X$ to equations of the form $f(x) = 0$, given a continuous function $f : X \to \RR$.
%			The set $\mathrm{dom}( \mathrm{Z}_X )$ of admissible instances of this problem is the set of all continuous functions $f$ with non-empty zero set $f^{-1}(0)$. The set $\mathrm{Z}_X(f)$ is the set of all zeros of $f$.
%		\end{example}

The benefit of using represented spaces is that mathematical problems can be simulated by functions acting on names:

\begin{defin}[Realizer]
Given represented spaces $(X, \delta), (Y, \gamma)$ and a problem $f : \subseteq X \rightrightarrows Y$, a $(\delta, \gamma)$-\emph{realizer} (or, simply, a \emph{realizer} when $\delta$ and $\gamma$ are clear) of $f$ is a function $F : \subseteq \baire \to \baire$ such that
    \[
\gamma (F(p)) \in f (\delta(p))
    \]
for every $p$ which is a $\delta$-name for an element of $\dom(f)$.
\end{defin}

\begin{defin}[Computable and realizer continuous problems]\label{Def_ComputableProblem}
    A problem $f : \subseteq (X, \delta) \rightrightarrows (Y, \gamma)$ is called \emph{$(\delta, \gamma)$-computable} (or simply \emph{computable}) if it has a computable realizer, and is called \emph{$(\delta, \gamma)$-realizer continuous} (or simply \emph{realizer continuous}) if it has a continuous $(\delta,\gamma)$-realizer.
\end{defin}

Since all computable functions on the Baire space are continuous, it follows that any computable problem between represented spaces is realizer continuous.
If $X$ and $Y$ are topological spaces we can ask whether the realizer continuous functions coincide with those that are topologically continuous.
The theory of representations provides an answer to this question using the notion of \emph{admissible representation}, usually employed for $T_0$ second countable represented spaces.
Admissible representations are those concretely used in computable analysis for such topological spaces. To be defined, this notion needs the introduction of a suitable relation among representations.

%We compare representations of sets by computable \emph{reductions}.

\begin{defin}
    Let $\delta_1, \delta_2$ be representations of a set $X$.
    A function $f : \subseteq \baire \to \baire$ \emph{reduces} (or \emph{translates}) $\delta_1$ to $\delta_2$ if $\delta_1(p) = \delta_2 (f(p))$ for all $p \in \dom(\delta_1)$.

    If $f$ is computable, we say that $\delta_1$ is \emph{computably reducible} to $\delta_2$ and write $\delta_1 \le \delta_2$. We write $\delta_1 \equiv \delta_2$ in case $\delta_1 \le \delta_2$ and $\delta_2 \le \delta_1$, and say that $\delta_1$ and $\delta_2$ are \emph{computably equivalent}.
    Similarly, we say that $\delta_1$ is \emph{continuously reducible} to $\delta_2$ when $f$ is continuous and write $\delta_1\le_c\delta_2$. The notation $\delta_1 \equiv_c \delta_2$ and the corresponding notion of being \emph{continuously equivalent} are defined accordingly.
\end{defin}

%Obviously, computable reducibility, resp.\ computable equivalence, implies continuous reducibility, resp.\ continuous equivalence.

\begin{defin}[Admissible representation]
Given a $T_0$ second countable space $X$, we say that $\delta: \subseteq \baire \to X$ is an \emph{admissible representation} for $X$ if it is continuous and all continuous representations of $X$ are continuously reducible to it.
\end{defin}

The fact that a representation $\delta_1$ is reducible to a representation $\delta_2$ means that $\delta_2$ does not contain more information than $\delta_1$; therefore, a representation is admissible if and only if it is, up to continuous equivalence, the \lq\lq poorest\rq\rq\ continuous representation of that space.

\begin{teo}[{\cite[Theorem 3.2.11]{weihrauch:2000}}]\label{Theorem_RealizerContinuous}
If $X$ and $Y$ are $T_0$ second countable spaces with admissible representations $\delta$ and $\gamma$, respectively, then $f: \subseteq  X \to Y$ is $(\delta, \gamma)$-realizer continuous if and only if it is topologically continuous.
\end{teo}

\subsubsection{Representation of computable metric spaces}\label{Subsection_RepresentationEuclideanSpace}

We view the Euclidean space $\RR^n$, with its dense subset $\QQ^n$ and the Euclidean metric, as an example of computable metric space:

\begin{defin}[Computable metric space]\label{Def_ComputableMetricSpaces}
    A \emph{computable metric space} $(X,d,\alpha)$ is a separable metric space $(X,d)$ with metric $d : X \times X \to \RR$ and a dense sequence $\alpha : \NN \to X$ such that $d \circ (\alpha \times \alpha) : \NN^2 \to \RR$ is a computable double sequence of real numbers, that is, such that $\langle n, m \rangle \mapsto d(\alpha_n, \alpha_m)$ is computable.
\end{defin}

In the case of the Euclidean space $\RR^n$, $\alpha:\NN\to\QQ^n$ is a canonical effective numbering of $\QQ^n$.

Computable metric spaces are usually represented by some version of Cauchy representation (see \cite[Lemma 4.1.6]{weihrauch:2000}), which is admissible for such spaces.

\begin{defin}
    The \emph{Cauchy representation} $\rho_X : \subseteq \baire \to X$ of a computable metric space $(X,d,\alpha)$ is defined by:
    \begin{align*}
        \dom(\rho_X) &:= \{\, p \in \baire \mid \lim_{n \to \infty} \alpha_{p(n)} \text{ exists and } \forall i\: d(\alpha_{p(i)} , \lim_{n \to \infty} \alpha_{p(n)} ) < 2^{-i}\, \} \\
        \rho_X(p) &:= \lim_{n \to \infty} \alpha_{p(n)}.
    \end{align*}
\end{defin}

%In the Example \ref{Example_Type2Machines}.\ref{Example_DecimalFractions} we showed a representation of real numbers for which multiplication by $3$ is not computable. This representation is not satisfactory for our purposes, as we want natural functions like arithmetic operations to be computable.
%We present now a representation of the real numbers which induces a reasonable notion of computability according to Definition \ref{Def_ComputableProblem}: the \emph{Cauchy representation} of the reals.

In particular, %we obtain the representation $\rho_\RR$ of real numbers, defined by:
%\begin{align*}
%    \rho_{\RR}( p ) = x \iff \forall i \in \NN \:\: |x - \alpha_{\QQ}(p(i))| < 2^{-i}.
%\end{align*}
%where $\alpha_{\QQ} : \NN \to \QQ$ is a canonical effective numbering of the rational numbers.
%In other words,
a Cauchy name of a real $x$ is an infinite sequence of rational numbers converging to $x$ at a fixed rate. We will sometimes identify the reals with their Cauchy names, so that by $x[j]$ we will denote the $j$-th rational number encoded in a given Cauchy name of $x$.
It is well known that all ordinary algebraic operations on the reals are computable with respect to $\rho_\RR$.

\subsubsection{Representation of continuous functions}

The set of all partial continuous functions $f : \subseteq \baire \to \baire$ has the same cardinality of the power set of $ \baire $, so it cannot be represented. This problem can be solved by considering only continuous functions with ``natural'' domains. By Kuratowski's Theorem (see \cite[Theorem 3.8]{kechris:1995}) for every partial continuous function $f : \subseteq \baire \to \baire$ there exists a $G_\delta$ set (i.e., a countable intersection of open sets) $G \subseteq \baire$ such that $\dom(f) \subseteq G$, and a continuous extension $g : G \to \baire$ of $f$. Therefore, it is not restrictive to consider only the set $\mathcal{C}$ of continuous functions with $G_\delta$ domain.

%\begin{defin}
    %We define the set $\mathcal{C}$ as the set of all continuous functions with $G_\delta$ domain, i.e.
   % \begin{align*}
     %   \mathcal{C} = \left\{ f : \subseteq \baire \to \baire \mid \mathrm{dom}(f) \text{ is a } G_\delta \text{ set} \right\}.
   % \end{align*}
We define then the following representation $\eta : \subseteq \baire \to \mathcal{C}$.
For $n \in \NN$, $p, q \in \baire$, let $\eta (\langle n, p \rangle) (q) $ denote the result computed by the Type-$2$ machine $M$ with code $n$ applied to the input $\langle p, q \rangle$.
We abbreviate $\eta (\langle n, p \rangle) $ with $\eta_{\langle n, p \rangle}$.
%\end{defin}
In fact, we may regard $M$ as an \emph{oracle machine}, which for every fixed oracle $p \in \NN^\NN$ computes the continuous function $\eta_{\langle n, p \rangle}$. Since every continuous function with $G_\delta$ domain can be computed by an oracle machine, $\eta$ really constitutes a representation for $\mathcal{C}$.	
%\begin{lemma}[\cite{weihrauch:2000} Lemma 2.3.11]
    %The function $\eta$ defined above is a representation of $\mathcal{C}$.
%\end{lemma}
Obviously, computable functions correspond exactly to those that have computable $\eta$-names.

Using the representation $\eta$ and Theorem \ref{Theorem_RealizerContinuous}, we obtain a representation for continuous functions between Euclidean (and actually much more general) spaces.
%\begin{lemma}[\cite{weihrauch:2000} Lemma 2.3.12]
   % A function $f : \subseteq \NN^\NN \to \NN^\NN$ is computable iff $f = \eta_{\langle x, p \rangle}$ for some $x \in \NN^*$ and some computable $p \in \NN^\NN$.
%\end{lemma}
Given a continuous function $f : \subseteq \RR^n \to \RR^m$, the idea is to use any $\eta$-name of a $(\rho_{\RR^n}, \rho_{\RR^m})$-realizer $F$ of $f$ as a name of $f$.
In particular, given $A \subseteq \RR^n$ we obtain a representation for
$$
C(A, \RR^m) = \{ f : \subseteq \RR^n \to \RR^m \mid f \text{ is continuous and } \dom(f) = A \}
$$
as follows (see \cite[Deﬁnition 6.1.1]{weihrauch:2000}):

\begin{defin}\label{Def_StandardRepresentationRealFunctions}
    The \emph{representation} $\delta^{\to}_A : \subseteq \baire \to C(A, \RR^m)$ is defined by
\[
        \delta^{\to}_A(p) = f \iff \forall q \left( \rho_{\RR^n} (q) \in A \implies f \circ \rho_{\RR^n}(q) = \rho_{\RR^m} \circ \eta_p(q) \right)
\]
    %As in Definition \ref{Def_RepresentedSpaceContinuousFunc}, when $A = \RR^m$, we simply write $\delta^{\to}$.
(i.e.\ $\eta_p$ is a realizer of $f$).
When $A = \RR^n$ and $n$ is clear from the context, we simply write $\delta^{\to}$.
\end{defin}

Since our partial functions have closed domains, we need to represent closed sets. There are three main representations used in Computable Analysis: negative, positive, and total. Intuitively, negative information consists in describing a closed set through a set of rational open balls covering its complement, the positive one in listing all rational open balls that intersect the closed set, and the total information is the combination of the previous ones.
As particular cases, when a closed set has a computable negative, positive, or total name, we obtain the co-r.e., r.e. or recursive closed sets, respectively. Hence, a closed r.e. set is a closed set for which the rational open balls intersecting it can be computably enumerated, and so on\footnote{For more details on representations and computability on closed and open sets of topological spaces one can see for example \cite{BP:2003}, \cite{Miller:thesis}.
An extensive comparison between different approaches to open sets in Reverse Mathematics is \cite{NS:2020}.}:

\begin{defin}[Closed set representations]\label{Def_ClosedSetRepresentations}
Let $(X, \delta, \alpha)$ be a computable metric space, and let $B_0,B_1,B_2....$ be a standard enumeration of its rational open balls, that is the open balls of the form $B(\alpha_n,r)$ for $n\in\NN$ and $r\in\QQ_+\cup\{0\}$.
By $\mathcal{A}_{-}(X)$ we denote the set of closed subsets of $X$ equipped with the \emph{negative information representation} $\psi_X^{-} : \subseteq \baire \to \{ F \subseteq X \mid F \text{ is closed in } X\}$ defined by
    \begin{align*}
        \psi_X^{-}(p) = F \iff F = X \setminus \bigcup_{i \in \ran(p)} B_i,
    \end{align*}
    where $\ran(p)$ denotes the range of $p: \NN \to \NN$.
    By $\mathcal{A}_{+}(X)$ we denote the set of closed subsets of $X$ equipped with the \emph{positive information representation} $\psi_X^{+} : \subseteq \baire \to \{ F \subseteq X \mid F \text{ is closed in } X\}$ defined by
    \begin{align*}
        \psi_X^{+}(p) = F 	&\iff \forall i \: \left( i+1 \in \ran(p) \iff B_i \cap F \ne \emptyset \right).
    \end{align*}
    Finally, by $\mathcal{A}(X)$ we denote the set of closed subsets of $X$ equipped with the \emph{total information representation} $\psi_X$ defined as $\psi_X:=\psi^+_X\wedge\psi^-_X$, that is
     \begin{align*}
        \psi_X (\langle p, q \rangle ) = F \iff \psi_X^{-}(p) =\psi_X^{+}(q)=F.
    \end{align*}
    %In other words, $\psi_X = \psi_X^{-} \wedge \psi_X^{+}$.
\end{defin}

The representations $\psi^-_X$, $\psi^+_X$ and $\psi_X$ are admissible with respect to the spaces $\mathcal{A}_-(X)$, $\mathcal{A}_+(X)$ and $\mathcal{A}(X)$ equipped with the upper-Fell, lower-Fell and Fell topology, respectively. It is clear from the definition of Cauchy representation that we can view $F$ as a computable element of $\mathcal{A}_{-}(X)$ if and only if we can ``semi-decide'' whether $x \not\in F$ for every $x \in X$.
This means that to show that a name for some $F \in \mathcal{A}_{-}(X)$ can be computed from some input $z$ it suffices to give a definition of $F$ by a $\Pi^0_1$-formula with parameter $z$.

When the dimension $n$ is clear from the context, we write $\psi^+$, $\psi^-$ and $\psi$ in place of $\psi^+_{\RR^n}$, $\psi^-_{\RR^n}$ and $\psi_{\RR^n}$, respectively.

The restriction of the representation $\psi^+$ to \emph{nonempty} closed sets $F \subseteq \RR^n$ is equivalent to the representation $\psi^+_e$ defined by
\begin{align}\label{Align_DenseSequencePositiveInformation}
    \psi_e^+ ( \langle p_0, p_1, \dots \rangle ) = F \iff \forall i \in \NN \: ( p_i \in \dom(\rho_{\RR^n}) ) \wedge F = \overline{ \{ \rho_{\RR^n}(p_i) \mid i \in \NN\} }.
\end{align}
%\todo[inline]{a volte usiamo $\&$ (con spaziature diverse) e a volte $\wedge$: sarebbe meglio uniformare}

In this paper we will usually encode closed sets by the total representation $\psi$.
A main feature of this representation is that it makes the distance function computable:

\begin{cor}\label{Cor_DistanceIsComputable}
The function $\Delta : \mathcal{A}(\RR^n) \setminus \{\emptyset\} \times \RR^n \to \RR_+$ defined by $\Delta(F,x) := d_F(x)=\inf\{d(x,y)\mid y\in F\}$ is $(\psi, \rho_{\RR^n}, \rho_\RR)$-computable.
\end{cor}

We are now ready to represent the set  $C_c(\RR^n)$ of all partial continuous functions from $\RR^n$ to $\RR$ with nonempty closed domain.
More precisely, $(f,F) \in C_c(\RR^n)$ if and only if $f$ is a partial continuous function over $\RR^n$ and $F = \dom(f) \in\mathcal{A}(\RR^n)\setminus\{\emptyset\}$:
\begin{align}\label{repparfunc}
\delta^\rightarrow_{C_c(\RR^n)}(\langle p,q\rangle) = f  \iff  \psi(q) = F \wedge \delta^\to_F(p) = f.
\end{align}

It is also useful to consider the following equivalent representation for $C_c(\RR^n)$, based on the fact that a continuous function (on a closed domain) is determined by all pairs of rational open balls $(B_i, B_j)$ such that $f(\overline{B_i}) \subseteq B_j$ (for $B_i$ intersecting its domain):

\begin{defin}
    %If $A \subseteq \RR^m$ is an open set, we define the representation $\delta'_A : \subseteq \baire \to C(A, \RR^n)$ by
    %\begin{align*}
        %\delta'_A (p) = f \iff \forall k = \langle i, j \rangle \in \NN \: \left( \, \overline{B_i} \subseteq A \land f(\overline{B_i}) \subseteq B_j \iff k \in \ran(p) \right).
   % \end{align*}
   %while
   %If $A \subseteq \RR^m$ is closed we define $\delta'_A$ by
    %\begin{align*}
       % \delta'_A (p) = f \iff \forall k = \langle i, j \rangle \in \NN \: \left( \, B_i\cap A \:\&\: f(\overline{B_i}) \subseteq B_j \iff k \in \ran(p) \right).
    %\end{align*}
       % \delta'_A (\langle p,q\rangle) = f \iff [\psi(q)=\mathrm{dom}(f)\:\&\\
   % \forall i\left( \, B_i \cap A \ne \emptyset \implies (f(\overline{B_i}) \subseteq B_j \iff k:=\langle i,j\rangle \in \ran(p)) \right) ].
$\delta'_{C_c(\RR^n)} (\langle p,q\rangle) = f$ if and only if $\psi(q) =\dom(f)$ and for all $k = \langle i, j \rangle$, $k \in \ran(p)$ exactly when $B_i\cap \dom(f)\neq\emptyset$ and $f(\overline{B_i}) \subseteq B_j$.
   % As in Definition \ref{Def_StandardRepresentationRealFunctions}, when $A = \RR^m$, we simply write $\delta'$.
\end{defin}

By naturally adapting \cite[Lemma 6.1.7]{weihrauch:2000}, one can indeed obtain

%\begin{lemma}[\cite{weihrauch:2000}, Lemma 6.1.7] \label{Lemma_EquivalenceRepresentationRealFunctions}
    %The representations $\delta^{\to}_A$ and $\delta'_A$ are equivalent, i.e.\ $\delta^{\to}_A \equiv \delta'_A$.
%\end{lemma}

\begin{lemma}\label{Lemma_EquivalenceRepresentationRealFunctions}
$\delta^\rightarrow_{C_c(\RR^n)}\equiv \delta'_{C_c(\RR^n)}$.
\end{lemma}

As a particular case, when only total continuous functions are considered, one can simplify  $\delta'_{C_c(\RR^n)}$ by omitting the components encoding the closed domains and obtain a representation equivalent to $\delta^\rightarrow$.

To deal with functions differentiable up to some order we introduce some useful shorthand.

\begin{notat}\label{Def_nInteger}
For $\bar{k} = (k_1, \dots, k_n) \in \NN^n$ let $|\bar{k}|:= k_1 + \dots + k_n \leq m$ and $\bar{k}! = k_1! \cdots k_n!$.
For $\bar{k} = (k_1, \dots, k_n)$ and $\bar{l}  = (l_1, \dots, l_n)$ we write $\bar{k} + \bar{l}$ for the $n$-tuple $(k_1 + l_1, \dots, k_n + l_n)$, and $\binom{\bar{k}}{\bar{l}}$ for the binomial coefficient, that is
\begin{align*}
\binom{\bar{k}}{\bar{l}} = \frac{\bar{k}!}{\bar{l}! (\bar{k} - \bar{l})!}.
\end{align*}
When we write $\bar{k} \le \bar{l}$ we mean that $\bar{k}$ is smaller than or equal to $\bar{l}$ component-wise.

Given $x \in \RR^n$ and $\bar{k} = (k_1, \dots, k_n)$, we write $x^{\bar{k}}$ for $x_1^{k_1} \cdots x_n^{k_n}$.
\end{notat}

The following easy observation will be useful.

\begin{prop}\label{Prop_(x-y)^l<=d(x,y)^l}
    For all $x, y \in \RR^n$ and $\bar{l} \in \NN^n$, $ |(x-y)^{\bar{l}} \,| \le d(x,y)^{|\bar{l}|}$.
\end{prop}
\begin{proof}
    It suffices to show that for every $x \in \RR^n$, $|x^{\bar{l}}| \le \lVert x \rVert ^{|\bar{l}|}$.
    This follows from $|x^{\bar{l}}|^2 = x_1^{2 l_1} \cdots x_n^{2 l_n}$ and $(\lVert x \rVert ^{|\bar{l}|})^2 =(\lVert x \rVert^2)^{|\bar{l}|}= (x_1^2 + \dots + x_n^2)^{l_1 + \dots + l_n}$.
\end{proof}

Now, for every $m > 0$, we introduce a representation for continuous and differentiable (total) real functions up to order $m$, i.e.,\ a representation for $C^m(\RR^n,\RR)$.
By Schwarz's Theorem a partial derivative of $f \in C^m(\RR^n,\RR)$ is uniquely determined by the sequence of natural numbers $\bar{k} := (k_1, \dots, k_n) \in \NN^n$ with $|\bar{k}| \leq m$ such that for every $1\leq i\leq n$ the element $k_i$ denotes the number of times the variable $x_i$ is derived (the order of derivation of the different variables being irrelevant).
Therefore, for $\bar{k} \in \NN^n$ we define $\partial_{\bar{k}}f :=\frac{\partial^{|\bar{k}|}f}{\partial x_1^{k_1}...x_n^{k_n}}$. We have then:
%Let $A \subseteq \RR^m$ be an open set,
%\begin{align*}
    %C^{k}(A) = \{ f : \subseteq \RR^m \to \RR \mid f \text{ is } C^{k} \text{ and } \mathrm{dom}(f) = A \}.
%\end{align*}
%Let $f \in C^{k}(A)$; since $A$ is open, the partial derivatives of $f$ are defined and continuous up to order $k$ in $A$.
%Therefore, we have the continuous functions
%\begin{align*}
    %f,\, \frac{\partial f}{ \partial x_1 }, \dots, \frac{\partial f}{ \partial x_m }, \, \frac{\partial^2 f}{ \partial x_1^2 },\, \frac{\partial^2 f}{ \partial x_1 x_2 },\, \dots, \frac{\partial^2 f}{ \partial x_1 x_m }, \dots , \frac{\partial^k f}{ \partial x_m^k }.
%\end{align*}
%In order to appropriately represent $f$ as a $C^{k}$ function, we need to provide information about the partial derivatives of $f$ up to order $k$.
%Let $K$ be the number of functions listed above, and consider the order in which they are listed; we can give the following definition.

%\begin{defin}\label{Def_StandardRepresentationCmFunctions}
   % The representation $\delta^{\to}_{C^{k}(A)} : \subseteq \baire \to C^{k}(A)$ is defined by
   % \begin{align*}
       % \delta^{\to}_{C^{k}(A)} \langle p_0, \dots, p_{K - 1} \rangle = f \iff \delta^{\to}_A (p_0) = f \& \dots \& \delta^{\to}_A (p_{K - 1}) = \frac{\partial^k f}{ \partial x_m^k }.
    % \end{align*}
    %When $A = \RR^m$, we simply write $\delta^{\to}_{C^{k}}$.
%\end{defin}

\begin{defin}\label{Def_StandardRepresentationCmFunctions}
The representation $\delta^{\to}_{C^m(\RR^n)} : \subseteq \baire \to C^m(\RR^n,\RR)$ is defined by
    \begin{align*}
\delta^{\to}_{C^m(\RR^n)} (\langle p_{\bar{k}} : |\bar{k}| \le m \rangle) = f \iff & \bigwedge_{|\bar{k}| \le m} \delta^\to( p_{\bar{k}} ) = \partial_{\bar{k}}f.
    \end{align*}
Here of course $\partial_{\bar{0}}f=f$ for $\bar{0}=(0,...,0)$.
\end{defin}
Sometimes we will use the notation $g^{(\bar{k})}$ to denote the partial derivative $\partial_{\bar{k}}g$, for some $\bar{k}$-derivable function $g$, to facilitate the comparison with jets and Whitney jets (Definitions \ref{Def_Jet} and \ref{Def_WhitneyJet}).

\section{Decomposition of open sets into cubes}\label{Section_DecompositionCubes}

In this section, we provide a computable construction of a family of cubes that covers the complement of a non empty closed set $F \subseteq \RR^n$, given with total information (Definition \ref{Def_ClosedSetRepresentations}); moreover the cubes of this family get smaller and smaller as they get closer to $F$.
The construction follows \cite{stein:1970}, but we need to introduce some nontrivial modifications in order to make every step computable.

Given $F \in \mathcal{A}(\RR^n) \setminus \{\emptyset\}$, we want to find a family $\FF = \left\{ Q_0, Q_1, ... \right\}$ of closed cubes with the following properties:
\begin{enumerate}[label=(C{\arabic*})]
    \item\label{Condition_C1} $\bigcup_{i \in \NN} Q_i = \RR^n \setminus
        F$,
    \item\label{Condition_C2} for all $i, j$ with $i \ne j$,
        $\textrm{int}(Q_i) \cap \textrm{int}(Q_j) = \emptyset$,
    \item\label{Condition_C3} there exist constants $c_1$ and $c_2$
        independent from $F$ such that for every $i$
        $$ c_1 \diam(Q_i) < d(Q_i, F) < c_2 \diam(Q_i). $$
\end{enumerate}

This will allow us, in the next section, to build a family of functions forming a ``partition of unity''.

First, we consider the set $\mathcal{M}_0$ of all cubes $Q \subseteq \RR^n$ whose edges have unit length and whose vertices have integer coordinates (obviously, the edges of these cubes are parallel to the coordinate axes).
Then, for all $k \ge 1$, we define $\mathcal{M}_k$ to be the set of all cubes obtained by splitting any $Q \in \mathcal{M}_{k-1}$ into $2^n$ new cubes via bisection of its edges, and define $\mathcal{M}_{-k}$ to be the set of cubes with edge length $2^{k}$ and whose vertex coordinates are multiple of $2^{k}$.
Observe that for all $k \in \ZZ$ and for every cube $Q \in \mathcal{M}_k$, the edges of $Q$ have length $2^{-k}$ and $\diam(Q) = \sqrt{n}\: 2^{-k}$.
We denote by $\mathcal{M}$ the union of the $\mathcal{M}_k$'s.

Now for any $Q \in \mathcal{M}$ and $i \in \NN$ let $R(Q,i)$ be the set of vertices of $Q$ together with the vertices of all cubes $Q' \subset Q$ that belong to $\mathcal{M}_{h_i}$, where $h_i$ is the least number that satisfies
$$ \frac{1}{2} \sqrt{n} \, 2^{-h_i} < 2^{-i}.$$
Notice that $R(Q,i) \subseteq Q$ is a finite set of points with rational coordinates such that for every $x \in Q$ there exists some $r \in R(Q,i)$ with $d(x,r) < 2^{-i}$.

\begin{remark}\label{Remark_R(Q,i)}
    Note that if $Q \in \mathcal{M}_k$, $Q' \in \mathcal{M}_{k+1}$ with $Q' \subset Q$, then for each $i, j \in \NN$, $R(Q,i)$ and $R(Q',j)$ have at least one point in common.
    In fact, since $Q'$ is obtained by splitting $Q$'s edges in half (notice that this also holds  for negative $k$), $Q$ and $Q'$ have at least one vertex in common.
\end{remark}

Each cube $Q \in \mathcal{M}$ can be represented by a pair $(i, c_Q)$ where $i \in \ZZ$ is such that $Q \in \mathcal{M}_i$ and $c_Q \in \QQ^n$ is the center of $Q$.
Fix now a computable bijection $\rho : \ZZ \times \QQ^n \to \NN$.

\begin{remark} \label{Remark_RepresentationCubes}
    Notice that, given  $i_1, i_2 \in \ZZ$, $y \in \QQ^n$ and $r \in \QQ$ we can uniformly compute the \emph{finite set} of all $Q \in \bigcup_{i_1 \le i \le i_2} \mathcal{M}_i$ with $c_Q \in B(y,r)$.
    Here and elsewhere by uniformly computing a finite set we mean that we can uniformly compute a natural number encoding it (in this case via $\rho$).
\end{remark}
		
%		By Remark \ref{Remark_RepresentationCubes} we can introduce an enumeration $\gamma: \NN \to \mathcal{M}$:
%		given a cube $Q \in \mathcal{M}$ with center $c_Q$ and sides of length $2^{-i}$ (i.e.\ $Q \in \mathcal{M}_i$), then $\gamma(n) = Q$ iff $\rho(i, c_Q) = n$.
		
We define a representation $\delta_\mathcal{M} : \baire \to \mathcal{P}(\mathcal{M})$ of $\mathcal{P}(\mathcal{M})$ by
%		\begin{align*}
%		\delta_\mathcal{M}(p) = \mathcal{F} \iff \mathcal{F} = \{ \gamma(n) \mid p(n) = 1 \}.
%		\end{align*}
 \begin{align*}
\delta_\mathcal{M}(p) = \bigcup_{i \in \ZZ} \{ Q \in \mathcal{M}_i \mid p(\rho(i,c_Q)) = 1 \}.
\end{align*}

With this representation, we can define the multi-valued function $\CUBE : \subseteq \mathcal{A}(\RR^n) \rightrightarrows \left( \mathcal{P}(\mathcal{M}), \delta_\mathcal{M} \right) $.
The admissible instances are all nonempty closed  $F \subseteq \RR^n$ given with total information, and $\CUBE(F)$ is the set of all  $\FF \in \mathcal{P}(\mathcal{M})$ that satisfy \ref{Condition_C1}, \ref{Condition_C2} and \ref{Condition_C3} when $c_1 = \frac12$ and $c_2=5$. We will show that CUBE is well defined and computable.

Next, we introduce two suitable families of cubes $\FF_0^F$ and $\FF^F$, and in order to make their construction uniformly computable we need to modify the definitions in \cite{stein:1970}. To this aim, for $k\in\ZZ$, we define two sequences of numbers:
%\begin{align}\label{Align_bk}
%    b_k := \begin{cases}
%        3k+1 &\text{ if } k > 0 \\
%        \mathrm{max}(k+3,0) &\text{ if } k \le 0
%    \end{cases} && \eta_k := 2^{-b_k}
%\end{align}
\begin{align}\label{Align_bk}
    b_k = \max(k+3,0), && \eta_k := 2^{-b_k}.
\end{align}
Notice that by the definition of the sequence $(b_k)_k$, for every $k \in \NN$, it holds that $b_k > k$.
This will be useful later on.

Let then
\begin{align*}
    \FF_0^F := & \bigcup_{k \in \ZZ} \left\{
    Q \in \mathcal{M}_k \middle\vert
    \begin{array}{l}
        \text{there exists } r \in R(Q, \max(k+1, 0)) \text{ satisfying}\\
        2 \diam(Q) - \eta_k < d(r,F)[b_k] < 4 \diam(Q) + \eta_k
    \end{array}
    \right\}.
\end{align*}
By Corollary \ref{Cor_DistanceIsComputable} the distance function $d(x,F)$ is computable. Moreover, $R(Q,i)$ is finite and uniformly computable. Finally, the inequalities defining $\FF_0^F$ can be reduced to inequalities between rational numbers. Hence, $\FF_0^F$ is \emph{uniformly computable} in $F$.

\begin{prop}\label{Prop_DistanceDiamRelation}
    For all $k \in \ZZ$ and $Q \in \mathcal{M}_k$, if $Q \in \FF_0^F$ then there is an $x \in Q$ such that
    \begin{align}\label{Align_Prop_DistanceDiamRelation}
        2 \diam(Q) - 2 \eta_k < d(x,F) < 4 \diam(Q) + 2 \eta_k.
    \end{align}
\end{prop}
\begin{proof}
    We can choose $x$ to be the $r$ witnessing $Q \in \FF_0^F$.
    Then the inequalities require the computation of an approximation of $d(x,F)$ up to $\eta_k = 2^{-b_k}$. By taking into account this round-off error we obtain the desired inequalities.
\end{proof}

\begin{prop}\label{Fact_EffectiveF}
    For each cube $Q \in \FF_0^F$, the set of cubes $Q' \in \mathcal{M}$ satisfying $Q \subseteq Q'$ and $Q' \in \FF_0^F$ is finite and uniformly computable.
\end{prop}
\begin{proof}
    Let $k$ be such that $Q \in \mathcal{M}_k$.
    By Proposition \ref{Prop_DistanceDiamRelation} there exists $y \in Q$ such that $d(y,F) < 4 \diam(Q) + 2 \eta_k$.
    If $Q'$ is a cube in $\mathcal{M}_h$ containing $Q$, with $h < k$ satisfying
    \begin{align}\label{Align_FactEffectiveF}
        \frac{1}{2} \diam(Q') = \frac{1}{2}\sqrt{n}\, 2^{-h} > 4 \diam(Q) + 2\eta_k,
    \end{align}
    then $Q'$ cannot belong to $\FF_0^F$.
    In fact, for each $x \in Q'$ (in particular if $x \in R(Q', \max(h+1, 0))$) it holds
    \begin{align*}
        d(x,F)[b_h] \le d(y,F) + d(x,y) + \eta_h \le \frac{1}{2} \diam(Q') + \diam(Q') + \eta_h.
    \end{align*}
    Using $b_h \geq h+2$ we have $2\eta_h = 2^{-b_h + 1} \leq 2^{-h-1} \le \sqrt{n} \, 2^{-h-1}  = \frac{1}{2} \diam(Q')$, so we obtain
    \begin{align*}
        d(x,F)[b_h] \le \frac{3}{2} \diam(Q') + \eta_h \le 2 \diam(Q') - \eta_h.
    \end{align*}
    So $Q'$ does not satisfy the conditions to be in $\FF_0^F$.
    Let now $h^*_k$ be the largest number satisfying (\ref{Align_FactEffectiveF}).
    Therefore a cube $Q'$ containing $Q$ can stay in $\FF_0^F$ only if $Q' \in \mathcal{M}_h$, with $h^*_k < h \le k$.
    To conclude note that there are only a finite number of cubes containing $Q$ and belonging to $\bigcup_{h^*_k < h \le k} \mathcal{M}_h$, and they can be detected effectively by definition of $\mathcal{M}$.
\end{proof}

We finally define $\FF^F \subseteq \FF_0^F$ as
\begin{align} \label{Align_FF}
    \FF^F := \left\{ Q \in \FF_0^F \mid \forall Q' \in \mathcal{M} \:\left( Q \subset Q' \to Q' \notin \FF_0^F \right) \right\}.
\end{align}
In other words $\FF^F$ consists of the maximal cubes of $\FF_0^F$. Each cube in $\FF_0^F$ has in fact a unique maximal cube in $\FF_0^F$ containing it.
Clearly these maximal cubes have disjoint interiors.

By virtue of Proposition \ref{Fact_EffectiveF} we can say that $\FF^F$ is uniformly computable from a name for $F$. In fact, in (\ref{Align_FF}), we can replace $\forall Q' \in \mathcal{M}$ with $\forall Q' \in \bigcup_{h^*_k < h \le k} \mathcal{M}_h$, where $k$ is such that $Q \in \mathcal{M}_k$. The set of such cubes $Q'$ containing $Q$ is finite, and uniformly computable.
		
Therefore, if we show that $\FF^F \in \CUBE(F)$, we have a computable realization for the multi-valued function $\CUBE$.

%In fact, given a $\psi$-name for $F$, thanks to Lemma \ref{Lemma_PsidEquivalentPsi}, we can computably find a name for $d_F$. %Using the name for $d_F$ we can find a $\delta_{\mathcal{M}}$-name for $\mathcal{F}_0^F$, since every operation in the %construction is computable.
In fact, given a $\psi$-name for $F$, we can computably find a $\delta_{\mathcal{M}}$-name for $\FF_0^F$, since every operation in the construction is computable (using Corollary \ref{Cor_DistanceIsComputable} and the definition of the finite sets $R(Q,i)$).
Then, using Proposition \ref{Fact_EffectiveF} we obtain a $\delta_{\mathcal{M}}$-name for $\FF^{F}$.
Therefore, by Definition \ref{Def_ComputableProblem}, $\CUBE$ is computable.

Since the set $F$ will always be clear from the context, from now on we write
$\FF$ and $\FF_0$ instead of $\FF^F$ and $\FF^F_0$.

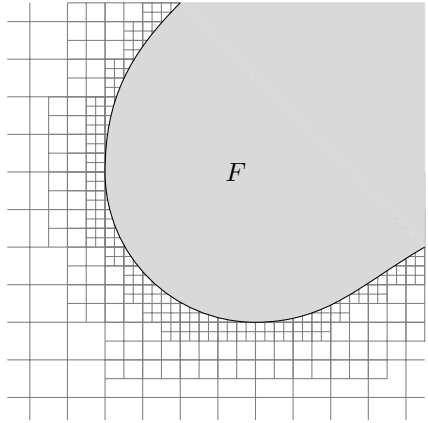
\begin{figure}[ht]
\centering
\begin{tikzpicture}[scale=1]
\draw[step=0.5cm,color=gray] (-0.3,-0.3) grid (5.25,5.25);

\draw[step=0.25cm,color=gray] (0.5,1) grid (2,5.25);
\draw[step=0.25cm,color=gray] (0.25,2) grid (1,4);
\draw[step=0.25cm,color=gray] (1,0.25) grid (4.75,3);
\draw[step=0.25cm,color=gray] (4.75,0.74) grid (5.25,3);

\draw[step=0.125cm,color=gray] (0.75,4) grid (1.25,2);
\draw[step=0.125cm,color=gray] (1,4.5) grid (1.5,4);
\draw[step=0.125cm,color=gray] (1.25,5) grid (1.75,4);
\draw[step=0.125cm,color=gray] (1.5,5.25) grid (2,5);
\draw[step=0.125cm,color=gray] (1,2) grid (1.75,1.75);
\draw[step=0.125cm,color=gray] (1.25,1.75) grid (2,1.25);
\draw[step=0.125cm,color=gray] (1.5,1.25) grid (4.25,1);
\draw[step=0.125cm,color=gray] (1.75,1) grid (4,0.75);
\draw[step=0.125cm,color=gray] (3,1.25) grid (4.75,1.75);
\draw[step=0.125cm,color=gray] (4.75,1.5) grid (5.25,2.5);

\draw[fill=gray!30, line width=0.1mm, gray!30] (2,5.25)
  -- (5.25,5.25)
  -- (5.25,2);
\draw[fill=gray!30] (2,5.25)
 	to [out=225,in=90] (1,3)
    to [out=270,in=180] (3,1)
    to [out=0,in=210] (5.25,2);

\draw (3,3) node[anchor=east]{$F$};

\end{tikzpicture}
\caption{The covering of cubes}\label{fig:Figure_coveringF}
\end{figure}
		
  It remains then to show that $\FF$ belongs to $\CUBE(F)$.
  As already seen, \ref{Condition_C2} is obvious from the construction.
  Since our definition of $\FF$ uses rational approximations and dense sets we need to prove \ref{Condition_C1} explicitly, while Stein can take it for granted.
		
\begin{prop}\label{Prop_FisCovering}
    The family of cubes $\FF$ covers the  complement of $F$, that is
    $$ \bigcup_{Q \in \FF} Q = \RR^n \setminus F. $$
\end{prop}
\begin{proof}

    If we prove that each point in $\RR^n \setminus F$ is contained in a cube in the family $\FF_0$, we are done.
    In fact, by definition of $\FF$, $\bigcup \FF = \bigcup \FF_0$, and by definition of $\FF_0$, no cube in $\FF_0$ intersects $F$.

    For each $k \in \ZZ$, recalling that $\sqrt{n} \, 2^{-k}$ is the diameter of the cubes in $\mathcal{M}_k$, we define the sets
    \begin{align*}
        C_k &= \{ x \notin F \mid 2 \sqrt{n} \, 2^{-k} \le d(x,F) < 2 \sqrt{n} \, 2^{-k} + 2^{-k-1} \}, \\
        D_k &= \{ x \notin F \mid 2 \sqrt{n} \, 2^{-k} + 2^{-k-1} \le d(x,F) \le 4 \sqrt{n} \, 2^{-k} - 2^{-k-1} \}, \\
        E_k &= \{ x \notin F \mid 4 \sqrt{n} \, 2^{-k} - 2^{-k-1} < d(x,F) < 4 \sqrt{n} \, 2^{-k} \}.
    \end{align*}
    Obviously, we have $\bigcup_{k \in \ZZ} \left( C_k \cup D_k \cup E_k \right) = \RR^n \setminus F$.
    %and also each $x \notin F$ is contained in exactly one of $C_k, D_k$ and $E_k$ for a single $k \in \ZZ$; in other words, $\bigcup_{k \in \ZZ} \{ C_k, D_k, E_k\}$ forms a partition of $\RR^n \setminus F$.

    Given $x \notin F$ there exists $k$ such that $x \in C_k \cup D_k \cup E_k$. Let $Q \in \mathcal{M}_k$ be a cube that contains $x$.

    If $x \in D_k$ let $r \in R(Q,\max(k+1,0))$ with $d(x,r) < 2^{-k-1}$. We then have
    \begin{align*}
        d(r,F)[b_k] \ge d(x,F) - d(x,r) - \eta_k > (2 \diam(Q) + 2^{-k-1}) - 2^{-k-1} - \eta_k, \\
        d(r,F)[b_k] \le d(x,F) + d(x,r) + \eta_k < (4 \diam(Q) - 2^{-k-1}) + 2^{-k-1} + \eta_k,
    \end{align*}
    so that $Q \in \FF_0$.\smallskip

    If $x \in C_k$ we consider two cases:
    \begin{itemize}
        \item[1)] There exists $r \in R(Q,\max(k+1,0))$ such that $d(r,F)[b_k] > 2 \diam(Q) - \eta_k$.
        \item[2)] Each $r \in R(Q,\max(k+1,0))$ satisfies $d(r,F)[b_k] \le 2 \diam(Q) - \eta_k$.
    \end{itemize}

    In the first case $Q$ is in $\FF_0$, in fact
    \begin{align*}
        d(r,F)[b_k] &< d(x,F) + d(x,r) + \eta_k  \\
        &< (2 \diam(Q) + 2^{-k-1}) + \diam(Q) + \eta_k < 4 \diam(Q) + \eta_k,
    \end{align*}
    since $2^{-k-1} < \sqrt{n} \, 2^{-k} = \diam(Q)$.

    In the second case, $Q$ is not in $\FF_0$, so we have to look for a smaller cube $Q'$. Let $Q'$ be a cube with $x \in Q'$, $Q' \in \mathcal{M}_{k+1}$ (and thus $2 \diam(Q') = \diam(Q)$), and $Q' \subseteq Q$.
    By Remark \ref{Remark_R(Q,i)} we have $R(Q,\max(k+1,0)) \cap R(Q',\max(k+2,0)) \ne \emptyset$.
    So, there exists $r \in R(Q', \max(k+2,0))$ such that
    \begin{align*}
        d(r,F)[b_{k+1}] &< d(r,F)[b_k] + \eta_{k}+\eta_{k+1} \le (2 \diam(Q) - \eta_k) + \eta_{k}+\eta_{k+1} \\
        &< 4 \diam(Q') + \eta_{k+1}.
    \end{align*}

    For the other inequality, by $x \in Q' \cap C_k$, it holds $d(r,F) \ge d(x,F) - d(x,r) \ge 2 \diam(Q) - \diam(Q') = 3 \diam(Q')$, so
    \begin{align*}
        d(r,F)[b_{k+1}] \ge 3 \diam(Q') - \eta_{k+1} > 2 \diam(Q') - \eta_{k+1}.
    \end{align*}	
    Therefore $Q' \in \FF_0$.\smallskip

    If $x \in E_k$ the situation is symmetric and the two cases are:
    \begin{itemize}
        \item[1)] There exists $r \in R(Q,\max(k+1,0))$ such that $d(r,F)[b_k] < 4 \diam(Q) + \eta_k$.
        \item[2)] Each point $r \in R(Q,\max(k+1,0))$ satisfies $d(r,F)[b_k] \ge 4 \diam(Q) + \eta_k$.
    \end{itemize}

    In the first case $Q \in \FF_0$, in fact
    \begin{align*}
        d(r,F)[b_k] &> d(x,F) - d(x,r) - \eta_k  \\
            &> (4 \diam(Q) - 2^{-k-1}) - \diam(Q) - \eta_k > 2 \diam(Q) - \eta_k.
    \end{align*}
			
    For the second case, let $Q'$ in $\mathcal{M}_{k-1}$, with $x \in Q'$ and $Q' \supseteq Q$. Then $\diam(Q') = 2 \diam(Q)$, and by Remark \ref{Remark_R(Q,i)} $R(Q,\max(k+1,0)) \cap R(Q',\max(k,0)) \ne \emptyset$.
    So there exists $r \in R(Q', \max(k,0))$ with
      \begin{align*}d(r,F)[b_{k-1}] &> d(r,F)[b_k] - \eta_k-\eta_{k-1} \ge (4 \diam(Q) + \eta_k) -\eta_k- \eta_{k-1}  \\
    &= 2 \diam(Q')- \eta_{k-1}.
    \end{align*}

    For the other inequality, by $x \in Q' \cap E_k$, it holds $d(r,F) \le d(x,F) + d(x,r) < 4 \diam(Q) + \diam(Q') = 3 \diam(Q')$, so
    \begin{align*}
        d(r,F)[b_{k-1}] < 3 \diam(Q') + \eta_{k-1} < 4 \diam(Q') + \eta_{k-1}.
    \end{align*}	
Thus $Q' \in \FF_0$ also in this case.
\end{proof}

To show that the cubes of $\FF$ gets smaller and smaller as they get closer to $F$ we show that condition \ref{Condition_C3} holds with $c_1 = \frac12$ and $c_2=5$.

\begin{prop}\label{Prop_DistanceDiamRelation2}
    For all $Q \in \FF$ it holds
    $$\frac{1}{2} \diam(Q) < d(Q,F) < 5 \diam(Q).$$
\end{prop}
\begin{proof}
    Let $k$ be such that $Q \in \mathcal{M}_k$.
    Consider now a witness $x \in Q$ of Proposition \ref{Prop_DistanceDiamRelation}, and $x' \in Q$ and $y \in F$ such that $d(Q,F) = d(x',y)$ ($x'$ and $y$ exist because $Q$ and $F$ are closed).
    We deduce that
    \begin{align*}
        d(Q,F) &= d(x',y) \ge d(x,y) - d(x,x') >  2\diam(Q) - 2\eta_k - \diam(Q)  \\
               &= \diam(Q) - 2\eta_k.
    \end{align*}
    Using $b_k > k+2$ (see (\ref{Align_bk})), we have $2\eta_k = 2^{-b_k + 1} < 2^{-k-1} \le \sqrt{n} \, 2^{-k-1}  = \frac{1}{2} \diam(Q)$, and thus $\frac{1}{2} \diam(Q) < d(Q,F)$.

    For the other inequality, let again $x \in Q$ be a witness to Proposition \ref{Prop_DistanceDiamRelation}, so that $d(x,F) < 4 \diam(Q) + 2\eta_k$.
    Using again $2\eta_k < \frac{1}{2} \diam(Q)$ we than have
    $$ d(Q,F) \le d(x, F) <  4 \diam(Q) + 2\eta_k < 5 \diam(Q).\qedhere$$
\end{proof}

We have thus completed the proof of:

\begin{teo}\label{Theorem_CUBEcomp}
The multi-valued function $\CUBE$ is computable.
\end{teo}

We now fix $\FF \in \CUBE (F)$.
We first obtain some further properties of $\FF$ and then define a derived family of enlarged cubes $\FF^*$.
To ease comparison with \cite{stein:1970} we adopt the following terminology.
		
\begin{defin}
    We say that two cubes \emph{touch} if their boundaries intersect.
\end{defin}

Since the cubes in $\FF$ have disjoint interiors by \ref{Condition_C2}, two cubes in $\FF$ intersect if and only if they touch.

\begin{prop} \label{Prop_DiamTouch}
    Let $Q_1 \in \mathcal{M}_h \cap \FF$ and $Q_2 \in \mathcal{M}_k \cap \FF$.
    If $Q_1, Q_2$ touch, then $|h-k| \le 3$, which is equivalent to
    \begin{align}\label{Align_DiamTouch8}
        \frac{1}{8} \diam(Q_1) \le \diam(Q_2) \le 8 \diam(Q_1).
    \end{align}
\end{prop}
\begin{proof}
    We show that if two cubes $Q_1,Q_2 \in \FF$ touch, then
    \begin{align}\label{Align_DiamTouch12}
        \frac{1}{12} \diam(Q_1) < \diam(Q_2) < 12 \diam(Q_1).
    \end{align}
    In fact, (\ref{Align_DiamTouch8}) follows immediately from (\ref{Align_DiamTouch12}) since the ratio of the diameters of two cubes in $\FF$ is always a power of $2$.

    Clearly it suffices to prove the second inequality in (\ref{Align_DiamTouch12}).
    By Proposition \ref{Prop_DistanceDiamRelation2} $\diam(Q_2) < 2 \, d(Q_2, F)$ and hence it is enough to prove $d(Q_2, F) < 6 \diam(Q_1)$.
    We achieve this goal by using the inequality $d(Q_1,F) < 5 \diam(Q_1)$ that holds by Proposition \ref{Prop_DistanceDiamRelation2}.
    In fact $d(Q_2,F) \leq d(x,F)$ for every $x \in Q_2$, in particular for every $x \in Q_1 \cap Q_2$.
    Fix then such an $x$, which exists by hypothesis, and let $x' \in Q_1$ and $y \in F$ be such that $d(x',y) = d(Q_1,F)$. Since $d(x,x') \leq \diam(Q_1)$ we have
    $$ d(Q_2,F) \leq d(x,y) \leq d(x,x') + d(x',y) < \diam(Q_1) + 5 \diam(Q_1) = 6 \diam(Q_1).\qedhere$$
\end{proof}

%\begin{prop} \label{Prop_BoundNumberTouches}
%    For every cube $Q \in \FF$ there are at most $24^n$ cubes in $\FF$ that touch $Q$, where $n$ is the dimension of the space.
%\end{prop}
%\begin{proof}
%    Let $Q \in \mathcal{M}_k$ for a given $k \in \ZZ$ and let $\mathcal{T}$ be the sets of all cubes in $\FF$ that touch $Q$.
%    By Proposition \ref{Prop_DiamTouch} $\mathcal{T} \subseteq \bigcup_{k-3 \leq i \leq k+3} \mathcal{M}_i$.
%
%    The cubes in $\mathcal{M}_{k+3}$ can be obtained by regularly decomposing all cubes of $\mathcal{M}_k$ into $8$-times rescaled copies of them. Any such decomposition results then in the production of $8^n$ new cubes for every decomposed cube in $\mathcal{M}_k$. Since there exist $3^n$ cubes in $\mathcal{M}_k$ touching $Q$ (including $Q$ itself), $Q$ touches at most $24^n$ cubes of $\mathcal{M}_{k+3}$.
%    Therefore, if $\mathcal T\subseteq\mathcal{M}_{k+3}$, then $\mathcal{T}$ contains at most $24^n$ elements. Suppose now that $\mathcal{T} \not\subseteq\mathcal{M}_{k+3}$. It exists then at least one cube $Q'\in\mathcal{T} \cap \mathcal{M}_i$ with $k-3 \le i < k+3$. By definition of $\FF$, no cube $Q''\in\mathcal{M}$ satisfying $Q''\subset Q'$ can belong to $\FF$, hence to $\mathcal T$. But $Q'$ contains several cubes $Q''\in\mathcal{M}_{k+3}$, hence one single $Q'\in\mathcal{M}_i$ added to $\mathcal T$ corresponds to several $Q''\in \mathcal{M}_{k+3}$ removed from $\mathcal T$. This means that the upper bound of $24^n$ cannot be exceeded in any case.
%\end{proof}

Fix a rational number $\varepsilon < \frac{1}{5}$.
Given $Q \in \FF$ let $c_Q$ be the center of $Q$ and define the \emph{enlarged cube} $Q^*$ as the image of $Q$ under the map $x \mapsto x^* := (1 + \varepsilon)(x - c_Q) + c_Q$.
Notice that $Q^*$ is $Q$ dilated by $(1 + \varepsilon)$.

%%		Given any cube $Q \in \FF$, let $Q^*$ be the expanded cube  $(1 + \varepsilon) (Q - c_Q) + c_Q$ for a \emph{fixed} $\varepsilon < 2^{-10}$, where $c_Q$ is the center of the cube $Q$.
%		%The choice to take $\varepsilon < 2^{-10}$ should be clear later.
We call the family of enlarged cubes $\FF^*$, i.e.
\begin{align} \label{Align_F*}
    \FF^* := \left\{ Q^* \mid Q \in \FF \right\}.
\end{align}
		
It is easy to see the following fact.

\begin{prop}\label{Prop_RelationQQ*}
    Each point $x \in Q$ is mapped to a point $x^* \in Q^*$ such that $d(x,x^*) \leq \frac{\varepsilon}{2} \diam(Q)$.
    Moreover $\diam(Q^*) = (1 + \varepsilon)\diam(Q)$, and for any $y \in \RR^n$ we have $d(y, Q^*) \ge d(y, Q) - \frac{\varepsilon}{2} \diam(Q)$.
\end{prop}

\begin{prop}\label{Prop_NonOverlapping}
    Every cube $Q^* \in \FF^*$ (exactly as the original cube $Q$) does not intersect $F$, and in fact $d(Q^*, F) \ge \frac{1}{2}( 1 - \varepsilon )\diam(Q)$.
\end{prop}
\begin{proof}
    Let $x\in Q$ and $y \in F$ be such that $d(x^*,y) = d(Q^*,F)$.
    Using Propositions \ref{Prop_DistanceDiamRelation2} and \ref{Prop_RelationQQ*}, we deduce that
    $$\frac{1}{2} \diam(Q) < d(Q,F) \leq d(x,y) \leq d(x^*,y) + d(x,x^*) \leq  d(Q^*,F) + \frac{\varepsilon}{2} \diam(Q).$$
    Therefore $0 < \frac{1}{2}( 1 - \varepsilon ) \diam(Q) \leq d(Q^*,F)$.
\end{proof}

\begin{defin} \label{Def_Fx}
    For each $x \notin F$ let
    \begin{align*}
        \FF_x = \{ Q \in \FF \mid x \in Q^* \}.
    \end{align*}
\end{defin}
		
Since the condition $x \in Q^*$ is not decidable we cannot compute $\FF_x$.
The following Proposition overcomes this obstacle by replacing, in contrast with the classical proof, $\FF_x$ with a suitable superset.
		
\begin{prop}\label{Prop_Gx}
    For each $x \notin F$, we can uniformly compute a finite set $\mathcal{G}_x$ such that $\FF_x \subseteq \mathcal{G}_x \subseteq \FF$.
    Moreover, the cardinality of $\mathcal{G}_x$ has a bound independent of $x$, indeed $|\mathcal{G}_x| \le N_n := \lceil 197 \, \sqrt{n} \rceil^n$.
\end{prop}
\begin{proof}
    Let $x \notin F$ and $\delta := d(x,F)$.
    If $Q \in \FF_x$ we have $d(Q^*, F) \le \delta \le d(Q^*, F) + \diam(Q^*)$.
    Using that for every $Q \in \FF$, $d(Q^*,F) < d(Q,F)$, and Propositions \ref{Prop_RelationQQ*} and \ref{Prop_DistanceDiamRelation2} we have
    \begin{align*}
        \delta \le d(Q^*, F) + \diam(Q^*) < d(Q, F) + (1 + \varepsilon) \diam(Q) < (6 + \varepsilon) \diam(Q).
    \end{align*}
    On the other hand, using Proposition \ref{Prop_NonOverlapping}, we have
    \begin{align*}
        \delta \ge d(Q^*, F) \ge \frac{1 - \varepsilon}{2} \diam(Q).
    \end{align*}
    Combining these two results (and using $\varepsilon < \frac13$), we obtain
    \begin{align*}
        \frac{\delta}{7} < \frac{\delta}{6 + \varepsilon} < \diam(Q) \le \frac{2 \, \delta}{ 1 - \varepsilon} < 3 \, \delta.
    \end{align*}

    Now, since obviously $d(x, c_Q) \le \frac{1}{2} \diam(Q^*)$,
    %, where $c_Q$ is the center of $Q$.
    using again Proposition \ref{Prop_RelationQQ*} (and $\varepsilon < \frac15$) we have
    \begin{align*}
        d(x,c_Q) \le \frac{1 + \varepsilon}{2} \diam(Q) \leq \frac{(1 + \varepsilon) \, \delta}{ 1 - \varepsilon} < \frac{3 \, \delta }{2}.
    \end{align*}

    We now are ready to define the computable set $\mathcal{G}_x$.
    First, we start computing $\delta$ using the total information of $F$ (see Corollary \ref{Cor_DistanceIsComputable}) until we find some $i$ such that $\delta[i] \ge 2^{-i + 1}$, so that we are sure that $\delta \ge 2^{-i}$.
    This always happens because $F$ is closed and therefore $\delta > 0$.
    Then, recalling that $\FF$ is uniformly computable from $F$, we define $\mathcal{G}_x$ to be the set of all cubes $Q \in \FF$ such that
    \begin{enumerate}
        \item $\frac{1}{14} \delta[i] < \diam(Q) < 6 \, \delta[i]$,
        \item $d(x,c_Q)[i] < \frac{3}{2} \, \delta[i] + 2^{-i + 2}$.
    \end{enumerate}

    By Remark \ref{Remark_RepresentationCubes}, since we have an upper bound and a lower bound on the edge length of the cubes in $\mathcal{G}_x$, and also a bound on the distance between the center of the cubes and $x$, we have that $\mathcal{G}_x$ is computable and can be listed in finite time.

    We still need to show that $\FF_x \subseteq \mathcal{G}_x$.
    Let $Q \in \FF_x$, then
    \begin{gather*}
        \frac{\delta[i]}{14} < \frac{\delta + 2^{-i}}{14} \le \frac{\delta}{7} < \diam(Q) < 3 \, \delta < 3 \,( \delta[i] + 2^{-i} ) \le 6 \, \delta[i], \\
        d(x, c_Q)[i] < d(x,c_Q) + 2^{-i} < \frac{3}{2}\, \delta + 2^{-i} < \frac{3}{2} \delta[i] + 2^{-i+2},
    \end{gather*}
    and $Q$ belongs to $\mathcal{G}_x$.

    Now, for the cardinality of $\mathcal{G}_x$, note first that the distance of any point in $\bigcup \mathcal{G}_x$ from $x$ is less than $\max_{Q \in \mathcal{G}_x} d(x, c_Q) + \max_{Q \in \mathcal{G}_x} \frac{1}{2} \diam(Q) < \frac{3}{2} \, \delta[i] + 5 \cdot 2^{-i} + 3 \, \delta[i]$.
    If we consider the cube $R$ centered in $x$ whose edges have length $14 \, \delta[i]$, then obviously $R \supseteq \bigcup \mathcal{G}_x$.
    To cover $R$ with the smallest cubes that can belong to $\mathcal{G}_x$ (Figure \ref{Figure_Gx}), we need
    \begin{align*}
        \left\lceil \frac{14 \,\delta[i]}{ \frac{1}{14 \, \sqrt{n}} \, \delta[i] } + 1 \right\rceil ^n = \lceil 196 \, \sqrt{n} + 1 \, \rceil ^n \le \lceil 197 \, \sqrt{n} \rceil^n
    \end{align*}
    of them. Therefore we have an upper bound for $|\mathcal{G}_x|$, that is $|\mathcal{G}_x| \le \lceil 197 \, \sqrt{n} \rceil^n$.
\end{proof}
		
		\begin{figure}[ht]
\centering
\begin{tikzpicture}[scale=1.8]

    \draw[fill=gray!30] (-1.5,2)
 	to [out=310,in=90] (-0.15,0.56)
    to [out=270,in=40] (-1.5,-0.5);

    \draw[step=0.25cm,color=gray!20] (-0.25,-1) grid (2.75,2);
	
	\draw[line width=0.2mm, black] (1,0.5) -- (1,1.5) -- (2,1.5) -- (2,0.5) -- (1,0.5);
	\draw[line width=0.2mm, gray, dashed] (0.8,0.3) -- (0.8,1.7) -- (2.2,1.7) -- (2.2,0.3) -- (0.8,0.3);
	\filldraw[black] (1.5,1) circle (0.25pt);
	
	\draw[line width=0.2mm, black] (1,0.5) -- (1,-0.5) -- (2,-0.5) -- (2,0.5) -- (1,0.5);
	\draw[line width=0.2mm, gray, dashed] (0.8,0.7) -- (0.8,-0.7) -- (2.2,-0.7) -- (2.2,0.7) -- (0.8,0.7);
	\filldraw[black] (1.5,0) circle (0.25pt);	
	
	\draw[line width=0.2mm, black] (0.5,0.5) -- (0.5,1) -- (1,1) -- (1,0.5) -- (0.5,0.5);
	\draw[line width=0.2mm, gray, dashed] (0.4,0.4) -- (0.4,1.1) -- (1.1,1.1) -- (1.1,0.4) -- (0.4,0.4);
	\filldraw[black] (0.75,0.75) circle (0.25pt);	
	
	\draw[line width=0.2mm, black] (0.5,0) -- (0.5,0.5) -- (1,0.5) -- (1,0) -- (0.5,0);
	\draw[line width=0.2mm, gray, dashed] (0.4,-0.1) -- (0.4,0.6) -- (1.1,0.6) -- (1.1,-0.1) -- (0.4,-0.1);
	\filldraw[black] (0.75,0.25) circle (0.25pt);	
	
	\filldraw[black] (1.06,0.45) circle (0.5pt)  node[anchor=north west]{$x$};
	
	\draw[line width=0.2mm, gray] (-0.25,-1) -- (-0.25,2) -- (2.75,2) -- (2.75, -1) -- (-0.25, -1);
	
	\draw (-1.15,0.5) node[anchor=west]{$F$};
	
	\draw (0,-0.7) node[anchor=west]{$R$};

\end{tikzpicture}

\caption{ } \label{Figure_Gx}

\end{figure}

%\begin{prop}\label{Fact_d(Q*y)>d(Q,y)}
%
%\end{prop}
%\begin{proof}
%    For every $x \in Q$, by Proposition \ref{Prop_RelationQQ*}, we have $d(x,y) \le d(x^*, y) + d(x, x^*) \le d(x^*, y) + \frac{\varepsilon}%{2} \diam(Q)$. Now, since $x \mapsto x^*$ is a bijection, we have
%    $$ d(y, Q) = \min_{x \in Q} d(x,y) \le \min_{x \in Q} d(x^*, y) + \frac{\varepsilon}{2} \diam(Q) = d(y, Q^*) + \frac{\varepsilon}{2} %\diam(Q).\qedhere $$
%\end{proof}

\begin{prop}\label{Prop_FydFinito}
    Given $x \notin F$, let $\delta = d(x,F)/2$. Then, the set %$\FF_{x,\delta}
    $\widetilde{\FF}_x := \{ Q \in \FF \mid d(x, Q^*) < \delta \}$ is finite.
\end{prop}
\begin{proof}
    First, observe that if $Q \in \widetilde{\FF}_x$, then $d(Q^*, F) < 3  \delta$.
    Moreover, we have $d(Q^*, F) + \diam(Q^*) > \delta$.
    In fact, let $x^* \in Q^*$ and $z \in Q^* \cap B(x, \delta)$, then
    \begin{align*}
        d(x^*,F) &\ge d(x,F) - d(x^*,x) \ge d(x,F) - d(x^*,z) - d(z,x) >\\
            &> 2 \delta - \diam(Q^*) - \delta = \delta - \diam(Q^*),
    \end{align*}
    hence $d(x^*, F) + \diam(Q^*) > \delta$ for every $x^* \in Q^*$, and so $d(Q^*, F) + \diam(Q^*) > \delta$.

    By Proposition \ref{Prop_NonOverlapping},
    we have
    \begin{align*}
        \frac{1 - \varepsilon}{2} \diam(Q) \le d(Q^*, F) < 3  \delta,
    \end{align*}
    and since $d(Q^*,F) < d(Q,F)$, using Propositions \ref{Prop_DistanceDiamRelation2} and \ref{Prop_RelationQQ*} we obtain
    \begin{align*}
        {\delta} < d(Q^*, F) + \diam(Q^*) < d(Q, F) + (1 + \varepsilon) \diam(Q) < (6 + \varepsilon)\diam(Q).
    \end{align*}
    Therefore we have upper and lower bounds for the diameter of the cubes $Q$ that belong to $\widetilde{\FF}_x$, i.e.\
    \begin{align*}
        \frac{\delta }{ 6 + \varepsilon} < \diam(Q) < \frac{6  \delta }{ 1 - \varepsilon}.
    \end{align*}

    Now, using Proposition \ref{Prop_RelationQQ*} we have
    \begin{align*}
        \delta > d(x, Q^*) \ge d(x, Q) - \frac{\varepsilon}{2} \diam(Q) > d(x, Q) - \frac{3  \delta  \varepsilon}{ 1 - \varepsilon },
    \end{align*}
    and using this inequality we obtain
    \begin{align*}
        \widetilde{\FF}_x \subseteq  \left\{ Q \in \FF \mid d(x, Q) < \frac{1 + 2\varepsilon}{1-\varepsilon}\delta \right\} := \mathcal{U}.
    \end{align*}

    Let now $k$ be the largest integer such that if $Q \in \mathcal{M}_k$ then $\diam(Q) > \delta / (6 + \varepsilon)$, that is, $k$ is the largest integer such that $\sqrt{n} \, 2^{-k} > \delta / (6 + \varepsilon)$. Using the fact that the cubes in $\mathcal{M}_k$ intersect only along the boundaries, have edges of length $2^{-k}$, and are the smallest possible among those that can belong to $\mathcal{U}$, we have $|\mathcal{U}| \le L^n$, where
    $$ L = \left \lceil \left( 2 \frac{1 + 2\varepsilon}{1 - \varepsilon}\delta + 2\frac{6 \delta}{1-\varepsilon} \right) 2^k + 1 \right \rceil = \left\lceil 2^{k} \frac{13 + 2\varepsilon}{1 - \varepsilon} \delta + 1\right\rceil.$$
    In fact, $L$ represents the maximum number of cubes along each edge of the grid formed by cubes in $\mathcal{M}_k$ whose union covers $\bigcup \, \mathcal{U}$ (see Figure \ref{Figure_FydFinito}).
\end{proof}

\begin{figure}[ht!]
\centering
\begin{tikzpicture}[scale=1.4]

	\draw[dashed] (0,1.25)
 	to [out=180,in=90] (-1.25,0)
    to [out=270,in=180] (0,-1.25)
    to [out=0,in=270] (1.25,0)
    to [out=90,in=0] (0,1.25);

	\draw[fill=gray!30] (0,1)
 	to [out=180,in=90] (-1,0)
    to [out=270,in=180] (0,-1)
    to [out=0,in=270] (1,0)
    to [out=90,in=0] (0,1);

    \draw[fill=gray!30] (-3,1.5)
 	to [out=310,in=90] (-2.1,-0.1)
    to [out=270,in=40] (-3,-1.75);

	\draw[step=0.25cm,color=gray!50] (-2.25,-2.25) grid (2.25,2.25);
	
	\draw[line width=0.2mm, black] (1,0.5) -- (1,1.5) -- (2,1.5) -- (2,0.5) -- (1,0.5);
	\draw[line width=0.2mm, black, dashed] (0.85,0.35) -- (0.85,1.65) -- (2.15,1.65) -- (2.15,0.35) -- (0.85,0.35);

	\draw (-0.4,0) node[anchor=west]{$B(x, \delta)$};
	
	\draw (-2.75,0) node[anchor=west]{$F$};
	
	\draw [decorate,
	decoration = {calligraphic brace}] (2.4,2.25) --  (2.4,-2.25);
	
	\draw [decorate,
	decoration = {calligraphic brace, mirror}] (-2.25, -2.4) --  (2.25,-2.4);
	
	\draw (2.6,0) node[anchor=west]{$L$};
	\draw (0, -2.6) node[anchor=north]{$L$};

\end{tikzpicture}

\caption{ } \label{Figure_FydFinito}

\end{figure}

\section{Partition of unity}\label{Section_PartitionUnit}
		
In this section, using the families of cubes $\FF$ and $\FF^*$ (defined respectively in (\ref{Align_FF}) and (\ref{Align_F*})), we define a collection of $C^{\infty}$ computable functions which constitutes a partition of unity on $\RR^n \setminus F$.
Moreover all partial derivatives of these functions are computable.
		
Let $Q_0 \subseteq \RR^n$ be the cube with center $(0,...,0)$ and edges parallel to the axes of length $1$.
We need to define a computable function $\varphi_0 \in C^\infty(\RR^n)$ such that $0 \leq \varphi_0(x) \leq 1$, $\varphi_0(x) = 1$ if $x \in Q_0$ and $\varphi_0(x) = 0$ for $x \notin Q_0^*$, and such that we can uniformly compute all partial derivatives of $\varphi_0$.

We postpone the actual definition of $\varphi_0$ and first we show how we use it.
First, we adapt $\varphi_0$ to any cube $Q \in \FF$ by defining $\varphi_Q \in C^\infty(\RR^n)$ as
\begin{align} \label{Align_varphi}
    \varphi_Q(x) := \varphi_0\left(\frac{x-c_Q}{l_Q}\right),
\end{align}
where $c_Q$ is the center of $Q$ and $l_Q$ the length of its edges. Notice that $0 \leq \varphi_Q(x) \leq 1$, $\varphi_Q(x) = 1$ if $x \in Q$ and $\varphi_Q(x) = 0$ for $x \notin Q^*$.
		
We then define the partial function $\varphi^*_Q \in C(\RR^n \setminus F)$ as
\begin{align}\label{Align_varphistar}
    \varphi^*_Q(x):=\frac{\varphi_Q(x)}{\Phi(x)}
\end{align}
for every $x \notin F$, where $\Phi(x) := \sum_{Q \in \FF} \varphi_{Q}(x)$.

\begin{remark} \label{Cor_SumVarphi}
For every $Q \in \FF$ let $d_Q \in \RR$.
Then
    \begin{align*}
        \sum_{ Q \in \FF} d_Q \varphi_Q(x) = \sum_{ Q \in \FF_x} d_Q \varphi_Q(x) = \sum_{ Q \in \mathcal{G}_x} d_Q \varphi_Q(x),
    \end{align*}
where $\FF _x$ and $\mathcal{G}_x$ are as in Definition \ref{Def_Fx} and Proposition \ref{Prop_Gx}.
The same applies to $\varphi^*_Q$ instead of $\varphi_Q$.
In fact if $Q \not\in \FF_x $ then $\varphi_Q(x) = 0$ and consequently also $\varphi_Q^*(x) = 0$.
\end{remark}
		
\begin{prop} \label{Prop_varphiCinfinito}
    For all $Q \in \FF$, $\varphi_Q$ and $\varphi_Q^*$ are $C^{\infty}$ and computable. Moreover, for all $Q \in \FF$ we can uniformly compute the partial derivatives of $\varphi_Q^*$.
\end{prop}
\begin{proof}
Computability and smoothness of $\varphi_Q$ follow from the same properties of $\varphi_0$.
By Remark \ref{Cor_SumVarphi} $\Phi(x) = \sum_{Q \in \mathcal{G}_x} \varphi_{Q}(x) $, and $\mathcal{G}_x$ is finite and computable by Proposition \ref{Prop_Gx}.
Moreover $\Phi(x) >0$ by Proposition \ref{Prop_FisCovering}.
Therefore $\varphi_Q^*$ is a quotient of computable $C^{\infty}$ functions with nonzero denominator and so is computable and $C^{\infty}$.

For the second part, let $Q \in \FF$ and let $1 \le i \le n$.
By the definition of $\varphi_Q^*$ (\ref{Align_varphistar}), we have
    \begin{align*}
         \partial_i \, \varphi_Q^*(x) = \frac{ \partial_i \, \varphi_Q(x) \, \Phi(x) - \varphi_Q(x) \, \partial_i \, \Phi(x) }{ \Phi^2(x) },
    \end{align*}
and $\partial_i \, \Phi(x) = \sum_{Q \in \mathcal{G}_x} \partial _i \, \varphi_Q (x)$. By the definition of $\varphi_Q$ (\ref{Align_varphi}) we have
    %$$\frac{ \partial  }{\partial x_i} \varphi_Q(x) = \frac{1}{l_Q} \, \frac{ \partial \, \varphi_0 }{\partial x_i} \left(\frac{x - c_Q}{l_Q} \right).$$
    $$ \partial_i  \, \varphi_Q(x) = \frac{1}{l_Q} \, \partial_i \, \varphi_0 \left(\frac{x - c_Q}{l_Q} \right).$$
Then, since we can compute the partial derivatives of $\varphi_0$, %by Proposition \ref{Prop_varphi0Derivatives}
we can also compute the partial derivatives of $\varphi_Q^*$.
For higher order derivatives, we proceed in the same way.
\end{proof}

By definition we have:

\begin{prop}\label{Prop_varphi*Partition}
The functions $\varphi^*_Q$ constitute a partition of unity on $\RR^n \setminus F$, that is
$$\sum_{Q \in \FF} \varphi^*_Q(x) = 1 \text{ for every } x \notin F.$$
\end{prop}
		
We now illustrate how to construct $\varphi_0$.
Since $\varphi_0$ depends on the dimension of the space, sometimes we refer to this function as $\varphi_0^n$, where $n > 0$ is the dimension of the space.
Let $\lambda : \RR \to \RR$ be the function
\[
    \lambda(x) := \begin{cases}
        e^{-\frac{1}{x}} & \text{ if } x > 0\\
        0 & \text{ if } x \le 0.\\
    \end{cases}
\]	
Let $\mu$ be the function defined by 		
\[
    \mu(x) := \frac{\lambda(x)}{\lambda(x) + \lambda(1-x)}.
\]
	
Now, using $\varepsilon$ defined in the previous section, we define $\nu$ by
\[
    \nu(x) := \mu\left( \frac{2x - \left( -1-\varepsilon \right) }{ \varepsilon } \right) \mu\left( \frac{ \left( 1 + \varepsilon \right) - 2x }{ \varepsilon } \right).
\]
This is a function that is $0$ outside $[-\frac{1}{2}-\frac{\varepsilon}{2},\frac{1}{2} + \frac{\varepsilon}{2}]$ and $1$ inside $[-\frac{1}{2},\frac{1}{2}]$.
To obtain $\varphi_0^n$ we set
$$ \varphi_0^n(x_1, ..., x_n) := \nu(x_1) \cdot ... \cdot \nu(x_n).$$

\begin{prop} \label{Prop_varphi0}
    The functions $\lambda$, $\mu$ and $\nu$ are $C^\infty$ and so is $\varphi_0^n$ for every $n>0$.
\end{prop}
\begin{proof}
    By \cite[Proposition 2.2]{nestruev:2020}, $\lambda$ is $C^\infty$, and since $\mu$ has a strictly positive denominator everywhere on the real line, $\mu$ is also $C^{\infty}$.
    Since $\nu$ is obtained by composition and product of $C^{\infty}$ functions, it is also $C^{\infty}$, and the same holds for $\varphi_0^n$ for all $n > 0$.
\end{proof}

\begin{lemma}\label{Lemma_Derivativesh}
    For every $k$, we can compute rational numbers $H_k$ and $B_k$ such that for every $x \in \RR$,
    \begin{align}\label{Align_Derivativesh}
\left| \lambda^{(k)}(x) \right| \le H_k, \quad \left| \mu^{(k)}(x) \right| \le B_k, \quad \text{and} \quad \left| \nu^{(k)}(x) \right| \le B_{k} \, \left( \frac{2}{\varepsilon}\right)^{k}.
    \end{align}
\end{lemma}
\begin{proof}
It is easy to check that for $x>0$ we have
    \begin{align*}
        \lambda^{({k})}(x) = \frac{\lambda(x)}{x^{2k}} P_k(x),
    \end{align*}
where $P_0(x) = 1$ and $P_{k+1}(x) = (1-2k x) P_k(x) + x^2 P'_k(x)$, so that $P_k$ is a polynomial with integer coefficients that can be explicitly computed.
It is an exercise to see that $\lambda(x)/x^{2k}$ for $x>0$ has maximum value $(2k)^{2k} e^{-2k}< (2k)^{2k}$.
Moreover when $x \in (0,1)$ $|P_k(x)| \leq A_k$ where $A_k$ is the sum of the absolute values of its coefficients.
Therefore we can set $H_k := (2k)^{2k} A_k$.

Now, for the derivatives of $\mu$, we have
\begin{align*}
        \mu^{(k)}(x) = \frac{Q_k \left( \left\{ \lambda(x), \lambda(1-x), \dots, \lambda^{(k)}(x), \lambda^{(k)}(1-x) \right\} \right) }{(\lambda(x) + \lambda(1-x))^{k+1}}
\end{align*}
 where $Q_k$ is a polynomial in $\left\{ \lambda(x), \lambda(1-x), \dots, \lambda^{(k)}(x), \lambda^{(k)}(1-x) \right\}$ with integer coefficients that can be computed.
    %Let $l_k$ be the number of terms of $Q_k$, and for every $j < l_k$, $q_j \in \ZZ$ be the coefficient of the $j$-th term.
    %Both $l_k$ and $(q_j)_{j < l_k}$ can be effectively computed using symbolic calculation for derivatives.
Furthermore, we can compute an upper bound  for the absolute value of each term of the polynomial using the bound for the derivatives of $\lambda$.
Therefore, we can compute a bound $T_k$ for the numerator by taking the sum of the upper bounds for each term.
For the denominator, we have that for every $x \in (0,1)$, since $\lambda$ is increasing we have $\lambda(x) + \lambda(1-x) \ge \lambda(1/2) = e^{-2} \geq 1/8$ and therefore
$(\lambda(x) + \lambda(1-x))^{k+1} \ge (1/8)^{k+1}$.
Summing up, if we set $B_k := 8^{k+1} T_k$ then $|\mu^{(k)}(x)| \le B_k$.

If we let $a = 1/2 + \varepsilon/2$ we have
% \[
%     \nu^{(k)}(x) = \frac{1}{\varepsilon ^k} \sum_{i \le k} (-1)^i \binom{k}{i} \mu^{(k-i)} \left(\frac{x+a}{\varepsilon} \right) \mu^{(i)} \left(\frac{-x+a}{\varepsilon} \right).
% \]
\[
    \nu^{(k)}(x) = \left(\frac{2}{\varepsilon} \right) ^k \sum_{i \le k} (-1)^i \binom{k}{i} \mu^{(k-i)} \left(\frac{2(x+a)}{\varepsilon} \right) \mu^{(i)} \left(\frac{2(-x+a)}{\varepsilon} \right).
\]
Notice that, since the derivatives of $\mu$ vanish outside $(0,1)$ and $\frac{2}{\varepsilon}(x+a) \in (0,1)$ if and only if $-\frac12 - \frac{\varepsilon}{2} < x <-\frac12$ while $\frac2\varepsilon(-x+a) \in (0,1)$ if and only if $\frac12 < x < \frac12 + \frac{\varepsilon}{2}$, only the terms
\[
    t_1:=\mu^{(k)} \left(\frac{2(x+a)}{\varepsilon} \right) \mu \left(\frac{2(-x+a)}{\varepsilon} \right) \quad \text{and} \quad t_2:=\mu \left(\frac{2(x+a)}{\varepsilon} \right) \mu^{(k)} \left(\frac{2(-x+a)}{\varepsilon} \right)
\]
do not vanish.
For every $x$, $\mu^{(k)} \left(\frac{2(x+a)}\varepsilon \right) \neq 0$ only when $\mu \left(\frac{2(-x+a)}\varepsilon \right)=1$, and analogously for the other term.
Moreover, the terms $t_1$ and $t_2$ cannot be non-zero simultaneously.
Therefore, the bound for $\nu^{(k)}$ holds as well.
\end{proof}

\begin{prop}\label{Prop_fghDerivatives}
The functions $\lambda$, $\mu$ and $\nu$ are computably $C^\infty$, i.e.\ given $k$ we can uniformly find $\delta^{\to}$-names for $\lambda^{(k)}$, $\mu^{(k)}$ and $\nu^{(k)}$.
\end{prop}
\begin{proof}
We first show that given $k$ and $x$ we can compute $\lambda^{(k)}(x)$ uniformly in $k$.
The main point here is that $e^{-1/x}$ and its derivatives are not defined in $0$, so we have to be careful when we do not know whether $x\le 0$ or $x>0$.

To show that $\lambda^{(k)}$ is computable
%Notice that we cannot apply Lemma \ref{Lemma_JoinFunctionComputable} because $e^{-1/x}$ is not defined on $0$.
%Therefore, we have to show that given $x \in \RR$, we can compute a good approximation of $f(x)$ until we realize, possibly, that $x \le 0$ or $x > 0$; in these cases, the output is $0$ or $e^{-1/x}$ respectively.
%In other words, we show that $f$ is $\delta^{\to}$-computable by defining a machine that given a $\rho$-name for $x$, compute a $\rho$-name for $f(x)$ (i.e.\ a $(\rho,\rho)$-realization of $f$).
we define a machine that operates in stages: at each stage $i$, it outputs an approximation of $\lambda^{(k)}(x)$ within an error of $2^{-i}$.
First let $j$ be least such that $2^{-j+1}<2^{-i}/H_{k+1}$.
If $x[j] \ge 2^{-i}/H_{k+1} -2^{-j}$ (and so $x \ge 2^{-i}/H_{k+1} -2^{-j+1}>0$), we proceed calculating an approximation of $\lambda^{(k)}(x)$ using the explicit formula for the $k$-th derivative of $e^{-1/x}$, while if $x[j] < 2^{-i}/H_{k+1} -2^{-j}$ (and so $x < 2^{-i}/H_{k+1}$), we print $0$ on the output tape.
Notice in fact that in the latter case, by the Mean Value Theorem, $|\lambda^{(k)}(x)| \le H_{k+1} x \le 2^{-i}$.
%: this is obvious if $x \le 0$.
%If instead $x>0$ we have $|f^{(k)}(x) - |x[i]|| = f^{(k)}(x) - |x[i]| < x - |x[i]| < 2^{-i}$ when $e^{-1/x} > |x[i]|$ (here we use $0< e^{-1/x} < x$ for every $x > 0$), while $|f(x) - |x[i]|| < |x[i]| < 2^{-i}$ when $e^{-1/x} \le |x[i]|$.

%The proofs for $\mu$ and $\nu$ are similar using $B_{k+1}$ in place of $H_{k+1}$.

The derivatives of  $\mu$ and $\nu$ can be computed from those of $\lambda$ using the familiar calculus rules.
\end{proof}

\begin{prop}\label{Prop_varphi0Derivatives}
    For every $n$ the function $\varphi_0^n$ is computably $C^\infty$.
\end{prop}
\begin{proof}
This is immediate from the previous Proposition because any partial derivative of $\varphi_0^n$ is the product of derivatives of $\nu$.
\end{proof}

We need upper bounds for the absolute values of the derivatives of the functions $\varphi_Q^*$; we postpone the proof of the next Proposition to the Appendix.

%Let $f,g : \RR^n \to \RR$ be smooth functions. We define
%\begin{align*}
%\mathrm{deg}^* \left(  \frac{\partial^{\bar{k}} \,f }{\partial x^{\bar{k}}} \, \frac{\partial^{\bar{l}} \,g }{\partial x^{\bar{l}}} \right) := |\bar{k}| + |\bar{l}|.
%\end{align*}

\begin{prop}\label{Prop_Derivativesvarphi}
Let $Q \in \FF$ and $l_Q$ be the length of its edges. For every $\bar{k}$, let
\begin{align}\label{Align_DefBbark}
    B_{\bar{k}} := B_{k_1} \cdot \dots \cdot B_{k_n}.
\end{align}%\label{Align_NormVarphilQ}
(where the $B_{k_i}$'s are defined as in Lemma \ref{Lemma_Derivativesh}).
Then, for every $x \notin F$,
\begin{align}
    \left| {\partial_{\bar{k}}} \varphi_Q(x) \right| \le B_{\bar{k}} \left(\frac{2}{ \varepsilon \, l_Q} \right) ^{|\bar{k}| }
\end{align}
Moreover,
\begin{align}\label{Align_NormVarphi^*DiamQ}
    \left| {\partial_{\bar{k}}} \varphi^*_Q(x) \right| \le  B'_{\bar{k}} \left( \frac2{ \varepsilon \diam(Q) }\right)^{|\bar{k}|},
\end{align}
where $B'_{\bar{k}}$ can be computed from $\bar{k}$ and the dimension of the space $n$.
\end{prop}

% 		\begin{cor} \label{Cor_varphiCinfinito}
% 			For all $Q \in \FF$, $\varphi_Q$ and $\varphi_Q^*$ are $C^{\infty}$ and computable. Moreover, for all $Q \in \FF$ we can computably find the partial derivatives of $\varphi_Q^*$.
% 		\end{cor}
% 		\begin{proof}
% 			For $\varphi_Q$, it is obvious. For $\varphi_Q^*$ instead, observe that $\Phi(x) = \sum_{Q \in \mathcal{G}_x} \varphi_{Q}(x) $, and $\mathcal{G}_x$ is finite and computable by Proposition \ref{Prop_Gx}. Using this we obtain that $\varphi_Q^* := \frac{\varphi_Q} {\Phi}$ is computable, since it is quotient of computable functions.
% 			
% 			For the second part, let $Q \in \FF$ and let $1 \le i \le n$. By the definition of $\varphi_Q^*$ (\ref{Align_varphistar}), we have
%			\begin{align*}
%				\frac{ \partial }{\partial x_i} \varphi_Q^*(x) = \frac{ \frac{ \partial  }{\partial x_i} \varphi_Q(x) \, \Phi(x) - \varphi_Q(x) \, \frac{ \partial  }{\partial x_i} \Phi(x) }{ \Phi^2(x) },
%			\end{align*}
%			and $\frac{ \partial  }{\partial x_i} \Phi(x) = \sum_{Q \in \mathcal{G}_x} \frac{ \partial  }{\partial x_i} \varphi_Q (x)$. By the definition of $\varphi_Q$ (\ref{Align_varphi}) we have
%			$$\frac{ \partial  }{\partial x_i} \varphi_Q(x) = \frac{1}{l_Q} \, \frac{ \partial \, \varphi_0 }{\partial x_i} \left(\frac{x - c_Q}{l_Q} \right).$$
%			 Then, since we can computably find the partial derivatives of the function $\varphi_0$ by Proposition \ref{Prop_varphi0Derivatives}, we can also computably find the partial derivatives of $\varphi_Q^*$. For higher order derivatives, we proceed in the same way.
% 		\end{proof}

\section{First Computable Whitney Extension Theorem}\label{Section_FirstWET}
	
In this section we prove what we call the first computable Whitney Extension Theorem. This is a refinement of the Tietze Extension Theorem stating the existence of a \emph{linear} operator mapping a continuous real-valued function defined on a closed subset of $\RR^n$ to a continuous total extension.
We show that Stein's proof of this result can be effectivized so that the resulting operator is  computable.
This provides also a proof of the computability of Tietze theorem which is  alternative to that given by Klaus Weihrauch (\cite{WeiTietze}), which is based on the effectivization of the proof contained in classical topology books (\cite{engelking:1989}).
Therefore, it turns out that both classical proofs approximate algorithmic constructions that are clearly revealed once their effectivization is pursued.

%Let $C_p(\RR^n, \RR)$ be the set of all partial continuous functions from $\RR^n$ to $\RR$. We let then $C_c(\RR^n) \subseteq C_p(\RR^n, \RR) \times \mathcal{A}(X)$ be the set of all partial continuous functions with nonempty closed domain. More precisely, $(f,F) \in C_c(\RR^n)$ iff $f$ is a partial continuous function over $\RR^n$ and $F = \dom(f) \neq \emptyset$.\marginpar{move to introduction?} Notice that $C_p(\RR^n, \RR)$, for mere cardinality reasons is not representable, whereas $C_c(\RR^n)$ is represented in the following way:
%\begin{align}
   % \delta_{C_c(\RR^n)}(\langle p,q\rangle) = f  \iff  \delta^\to(p) = f \wedge \psi(q) = \dom(f).
%\end{align}
		
In the proof of the Whitney Extension Theorem given in \cite{stein:1970}, the author, after defining his own family of cubes and partition of unity, defines the operator $\mathcal{E}_0 $ which sends $f$ into a total continuous extension, as follows:
\begin{align}
    \mathcal{E}_0(f)(x) := \begin{cases}
                                f(x) &\textrm{ if }x \in F \\
                                \sum_{Q \in \FF} f(p_Q) \varphi^*_Q(x) &\textrm{ else}
                            \end{cases}
\end{align}
where, for every $Q \in \FF$, $p_Q$ is a \emph{projection point} of the cube $Q$ on $F$, i.e.\ a point in $F$ such that $d(Q,F) = d(p_Q, Q)$.
In Sections \ref{Section_DecompositionCubes} and \ref{Section_PartitionUnit} we gave an effective construction of our $\FF \in \CUBE(F)$ and of the partition of unity $\{ \varphi^*_Q \mid Q \in \FF\}$.

The next step to effectivize the proof is to replace projections with approximate projections.
First, we determine upper bounds for $d(x^*,p_Q)$ for every $Q \in \FF$ and $x^*\in Q^*$; this will help us substituting exact projections with approximate projections.
By Proposition \ref{Prop_NonOverlapping} we know that $\frac{1}{2} (1 - \varepsilon )\diam(Q) \leq d(Q^*,F)$.
We set $e := \frac {2}{1 - \varepsilon }$ so that $\diam(Q) \leq e \, d(Q^*,F).$ Using this, we prove the following
		
\begin{prop}\label{Fact_RelationDistanceProjection}

    For any cube $Q \in \FF$, $x^* \in Q^*$ and $y \in F$, %and projection point $p_Q$ on $F$,
    $d(x^*,p_Q) < 7 e \, d(x^*,y)$.
\end{prop}
\begin{proof}
%    Notice that $x^*$ is furthest from $p_Q$ when the closest point in $Q$ to $p_Q$ is a vertex $x'$, and $x^*$ is the vertex of $Q^*$ opposite to $x'$, that is, $x^*$ is the image of the vertex $x$ of $Q$ opposite to $x'$ (as shown in the following figure). 	
Let $x$ be the preimage of $x^*$ and $x'\in Q$ be a point such that $d(p_Q, x') = d(p_Q, Q)$; then, using {Propositions \ref{Prop_DistanceDiamRelation2} and \ref{Prop_RelationQQ*}}, we obtain
    \begin{align*}
        d(x^*,p_Q) &\leq d(x^*,x) + d(x,x') + d(x',p_Q) \leq \frac{\varepsilon}{2} \diam(Q) + \diam(Q) + d(Q,F) \\
        &< 7 \, \text{diam}(Q) \leq 7 e \, d(Q^*, F) \le 7 e \, d(x^*, y).\qedhere
    \end{align*}
\end{proof}

Since, as shown in \cite{proj:2019}, exact projections cannot in general be computed, we replace $p_Q$ with a point $r_Q \in F$ such that $d(r_Q,Q) < 5 \diam(Q)$, so that the upper bound just provided still holds.

\begin{lemma}\label{Lemma_RelationDistanceApproxProjection}
    For every cube $Q \in \FF$, we can compute a point $r_Q \in F$ with $d(r_Q,Q) < 5 \diam(Q)$. Moreover, for all $x^* \in Q^*$ and $y \in F$, $d(x^*,r_Q) < 7 e \, d(x^*,y)$.
\end{lemma}
\begin{proof}
Positive information gives us a dense sequence in $F$ (according to formula (\ref{Align_DenseSequencePositiveInformation})), while $\bigcup_{i \in \NN} R(Q,i)$ is a dense sequence in $Q$.
We search for $r_Q$ in the dense sequence in $F$ and $q \in \bigcup_{i \in \NN} R(Q,i)$ such that $d(r_Q, q) < 5 \diam(Q)$ (our search succeeds by Proposition \ref{Prop_DistanceDiamRelation2}).

As in the proof of the Proposition \ref{Fact_RelationDistanceProjection} we obtain $d(x^*,r_Q) \le 7 e \, d(x^*, y)$ using as $x'$ the closest point in $Q$ to $r_Q$.
\end{proof}

Consequently, we modify the definition of $\mathcal{E}_0$ using the points $r_Q$.
This modification leads to define the total continuous extension of $f$ by
\begin{align} \label{Align_Def_g(x)}
    g(x):=
    \begin{cases}
        f(x) 	& \textrm{if } x \in F \\
        \sum_{Q\in\FF} f(r_Q) \varphi^*_Q(x) & \textrm{else}.
    \end{cases}
\end{align}

Since both $\FF$ and the points $r_Q$ depend on the $\psi$-name of $F$, we actually obtain a  multi-valued function $\mathrm{WET}_0 : C_c(\RR^n) \rightrightarrows C(\RR^n)$ (recall that $C_c(\RR^n)$ is the space of continuous partial functions $f:\subseteq\RR^n\to\RR$ with non-empty closed domain).
If we fix a $\psi$-name $q$ for $F$, we obtain a single valued map $\mathrm{WET}^{\,q}_0$ which is linear.
In fact, $q$ determines $\FF$ (given by the computable realizer of $\CUBE$), the family of functions $\{ \varphi_Q^* \mid Q \in \FF \}$  and also the point $r_Q$ for all $Q \in \FF$.
Then, by definition (\ref{Align_Def_g(x)}), the operator $\mathrm{WET}^{\,q}_0$ is clearly single-valued and linear.\smallskip

To prove that $\mathrm{WET}_0$ is computable we need to know how much a closed ball $\overline{B(y',\delta)}$ containing a point $y \in F$ and a point $x \notin F$ should be enlarged in order to include $r_Q$ for every $Q \in \FF_x$.
To this purpose, using Lemma \ref{Lemma_RelationDistanceApproxProjection} and $x \in Q^*$, we obtain
\begin{equation}
\label{Align_d(y',rQ)}
    d(y',r_Q) \leq d(y',x) + d(x,r_Q) \leq \delta + 7e \, d(x,y) \leq \delta + 7e \, 2 \delta = (14e + 1)\delta.
\end{equation}
		
Recall that we represent $C_c(\RR^n)$ by the representation $\delta^\rightarrow_{C_c(\RR^n)}$ defined in formula (\ref{repparfunc}).  We are now ready to prove the first computable Whitney Extension Theorem:

\begin{teo}[First Computable Whitney Extension Theorem] \label{Theorem_FirstWET}
    The multi-valued function $\mathrm{WET}_0 : C_c(\RR^n) \rightrightarrows C(\RR^n)$ defined above is computable.
\end{teo}
\begin{proof}
We give an algorithm computing the $\mathrm{WET}_0$ extension $g$ for any given input made of a name for $f$ and a total description of $F = \dom(f) \neq \emptyset$.
%We inspect the total description of $F$, in particular the positive information, and if we find out that $F=\emptyset$  (that is, the first digit of the sequence encoding the positive information is $0$, see Remark \ref{Remark_PositiveInformationEmptyset}), then the constant function equal to $0$ is released as output. \marginpar{Verificare se togliere il caso $F = \emptyset$}
%Otherwise we describe an algorithm to compute $g(x)$ for every $x \in \RR^n$, which means that we have a method to find a standard name of the total extension $g$.
By Lemma \ref{Lemma_EquivalenceRepresentationRealFunctions}, we can assume that the name for $f$ encodes a list of pairs of basic open sets in $\RR^n$ satisfying the following requirement: $(J,L)$ occurs in the list if and only if $J \cap F \neq \emptyset$ and $f(\overline J)\subseteq L$.

Given $x \in \RR^n$ we need to compute better and better approximations of $g(x)$, that is, at  stage $i$ we must output a rational number $q$ satisfying $|q - g(x)|\leq 2^{-i}$.
To this end, on one hand we inspect the negative information of $F$ and if we find an open ball $B \subseteq \RR^n \setminus F$ with $x \in B$, then we compute $\sum_{Q \in \FF} f(r_Q)\varphi^*_Q(x)$, according to (\ref{Align_Def_g(x)}).	
As long as this search does not succeed, we try to apply the function $f$ to $x$, as if we believe $x \in F$.
If $x \in F$ really holds, then the application of $f$ to $x$ will give the correct result, but if $x \notin F$ we will eventually know it and suddenly realize that the correct output is $\sum_{Q \in \FF}f(r_Q)\varphi^*_Q(x)$.

Let $c := 14 e + 1$ and, at stage $i$, search in the information about $f$ for a pair of the form $(B(y',c \, \delta),B(y'',\epsilon))$ such that
    \begin{enumerate}
        %\item\label{Item_FirstWET_1} if $i > 0$, $B(y',\delta)\subseteq J_{i-1}$, $B(y'',\epsilon)\subseteq L_{i-1}$,
        \item\label{Item_FirstWET_2} $x \in B(y',\delta)$ and $B(y',\delta)$ is listed in the positive information of $F$,
        \item\label{Item_FirstWET_3} $\epsilon \leq 2^{-i}$.
    \end{enumerate}
Notice that $x \in F$ is a sufficient, although not necessary, condition for such a pair to exist.
If $(B(y',c \, \delta),B(y'',\epsilon))$ is in the information of $f$ (even if $x \notin F$) then by definition of the constant $c$ and by (\ref{Align_d(y',rQ)}) not only $x \in \overline{B(y',c \, \delta)}$, but also $r_Q \in \overline{B(y', c \, \delta)}$ for every $Q \in \FF_x$.
This means that $f(r_Q) \in B(y'', \epsilon)$ for all such cubes $Q$. 	

When we find $(B(y',c \, \delta), B(y'',\epsilon))$ satisfying the requirements, we write $y''$ on the output tape.
Now, we show that $y''$ provides a good approximation for $g(x)$, i.e.\ $|y'' - g(x)| < 2^{-i}$.

If $x \in F$, then $|y'' - f(x)| < \epsilon \leq 2^{-i}$, since $x \in B(y',\delta) \subseteq B(y',c \, \delta)$ and by assumption $f(\overline{B(y',c  \delta)}) \subseteq B(y'', \epsilon)$.		
If instead $x \notin F$, using Remark \ref{Cor_SumVarphi} and Proposition \ref{Prop_varphi*Partition} we deduce
\begin{align*}
    |y''-g(x)| &= \left| y''-\sum_{Q \in \FF} f(r_Q) \varphi^*_Q(x) \right| = \left| y''-\sum_{Q\in\FF_x}f(r_Q)\varphi^*_Q(x) \right|  \\
    &= \left| y'' \sum_{Q \in \FF_x} \varphi^*_Q(x) - \sum_{Q\in\FF_x} f(r_Q) \varphi^*_Q(x) \right|  \displaybreak[0]\\
    &= \left| \sum_{Q \in \FF_x}(y'' - f(r_Q)) \varphi^*_Q(x) \right| \leq \sum_{Q \in \FF_x} \left| (y'' - f(r_Q)) \varphi^*_Q(x) \right|  \displaybreak[0]\\
    & =\sum_{Q \in \FF_x} |y'' - f(r_Q)|\varphi^*_Q(x) < \sum_{Q \in \FF_x} \epsilon \, \varphi^*_Q(x) \\
    &= \epsilon \sum_{Q \in \FF_x} \varphi^*_Q(x)= \epsilon \leq 2^{-i}.
\end{align*}

If the pair $(B(y',\delta),B(y'',\epsilon))$ does not exist, then we will eventually realize that $x \notin F$, and when this happens we can compute an approximation of $g(x) = \sum_{Q \in \FF} f(r_Q)\varphi^*_Q(x)$.
\end{proof}

\section{Computable Whitney Extension Theorem}\label{Section_WETm}
		
In this section we want to extend Theorem \ref{Theorem_FirstWET} to the case of $C^m$ functions, for $m \in \NN$. First, we need to define what we mean by differentiable functions on a closed subset of $\RR^n$. To do this, we introduce the \emph{Whitney jets}, which are finite collections of functions that behave as partial derivatives of each other.
Recall this classical theorem from Real Analysis.

% \begin{teo}[Multivariate Taylor's Theorem with Peano Remainder] \label{Theorem_MultivariateTaylor}
%     Let $f : \RR^n \to \RR$ be a $m$-times differentiable function at the point $y \in \RR^n$. Then, there exists a function $R(x,y)$ which is $o( d(x,y)^m )$ and for all $x \in \RR^n$ satisfies
%     \begin{align}
%         f(x) = \sum_{|\bar{k}| \le m} \frac{1}{\bar{k}!} \, { \partial _ {\bar{k}}} f(y) \, (x - y)^{\bar{k}} + R(x,y).
%     \end{align}
% \end{teo}

\begin{teo}[Multivariate Taylor's Theorem with Lagrange Remainder] \label{Theorem_MultivariateTaylor}
    Let $f : \RR^n \to \RR$ be an $(m+1)$-times continuously differentiable function at the point $y \in \RR^n$. Then, there exists a function $R(x,y)$ which is $O( d(x,y)^{m+1} )$ and for all $x \in \RR^n$ satisfies
    \begin{align}
        f(x) = \sum_{|\bar{k}| \le m} \frac{1}{\bar{k}!} \, { \partial _ {\bar{k}}} f(y) \, (x - y)^{\bar{k}} + R(x,y).
    \end{align}
    % where
    % $$ R(x,y) = \sum_{|\bar{l}| = m+1} \frac{1}{\bar{l}!} \, { \partial _ {\bar{l}}} f(\xi) \, (x - y)^{\bar{l}} $$
    % and $\xi$ is a point in the segment $[x,y]$.
\end{teo}

We now introduce the notion of \emph{Whitney jet} over a closed subset of $\RR^n$.

\begin{defin}\label{Def_Jet}
    Let $F \subseteq \RR^n$ be a closed set, and $m \in \NN$. A \emph{jet} of order $m$ on $F$ is an indexed set of  functions $(f^{(\bar{k})})_{|\bar{k}| \le m}$ with domain $F$.
\end{defin}

Given a jet of order $m$ on $F$, $(f^{(\bar{k})})_{|\bar{k}| \le m}$, and $y \in F$, we define the functions $P^{\bar{k}}_y  : \RR^n \to \RR$ as
\begin{align}\label{Align_P(x,y)}
    P^{\bar{k}}_y(x) := \sum_{|\bar{k} + \bar{l}| \le m} \frac{1}{ \bar{l}! } {f^{(\bar{k} + \bar{l})} (y) (x - y)^{\bar{l}}} .
\end{align}
Note that, if $y$ is fixed, these functions are polynomials in $x$.
For ease of notation we use $P_y(x)$ for $P^{\bar{0}}_y(x)$.

\begin{defin}\label{Def_WhitneyJet}
    A \emph{Whitney jet} of order $m$ on $F$ is a jet $(f^{(\bar{k})})_{|\bar{k}| \le m}$ of order $m$ on $F$ where for each $\bar{k}$, $f^{(\bar{k})}$ is continuous on $F$, and there exists $M > 0$ such that for every $\bar{k}$
    \begin{equation}\label{Align_CompatibilityConditions}
    \begin{split}
        &f^{(\bar{k})}(x) = P^{\bar{k}}_y(x) + R^{\bar{k}}(x,y) \quad\quad \forall x, y \in F, \\
        &\text{where } R^{\bar{k}}(x,y) \text{ satisfies } |R^{\bar{k}}(x,y)| \le M \, d(x,y)^{m-|\bar{k}| + 1}.
    \end{split}
    \end{equation}
    For ease of notation we use $R(x,y)$ for $R^{\bar{0}}(x,y)$.

    We call the conditions in (\ref{Align_CompatibilityConditions}) \emph{compatibility conditions}.
    We write $\mathcal{J}^m_c(\RR^n)$ for the set of all Whitney jets of order $m$ with closed domain in $\RR^n$, and, given a closed set $F \subseteq \RR^n$, we write $\mathcal{J}^m(F)$ for the set of all Whitney jets of order $m$ on $F$.
\end{defin}

\begin{defin}
    We now define the representation of a Whitney jet $(f^{(\bar{k})})_{|\bar{k}| \le m}$ of order $m$ on $F$ as
    \begin{align*}
        \delta_{\mathcal{J}^m_c(\RR^n)} (\langle \langle p_{\bar{k}} : |\bar{k}| \le m \rangle , q, M  \rangle) = (f^{(\bar{k})}) \iff & \bigwedge_{|\bar{k}| \le m} \delta^\to( p_{\bar{k}} ) = f^{(\bar{k})} \wedge \psi(q) = F \\ \wedge \: &(M\in \NN \text{ satisfies the conditions \ref{Align_CompatibilityConditions})}.
    \end{align*}
\end{defin}

\begin{prop} \label{Prop_DlPj=Pj+l}
    Given a Whitney jet $(f^{(\bar{k})})_{|\bar{k}| \le m}$ of order $m$ on $F$, $\bar{k}$ and $\bar{h}$ with $|\bar{k} + \bar{h}| \le m$, then for all $x \in \RR^n$ and $y \in F$, it holds
    \begin{align*}
        { \partial_{\bar{h}} } \, P^{\bar{k}}_y (x) = P^{ \bar{k} + \bar{h} }_y (x).
    \end{align*}
\end{prop}
\begin{proof}
We have
    \begin{align*}
        { \partial_{\bar{h}} } P^{\bar{k}}_y (x) &= \sum_{|\bar{k} + \bar{l}| \le m} \frac{f^{(\bar{k} + \bar{l})}(y)}{ \bar{l}! } \partial_{\bar{h}} (x - y)^{\bar{l}}  \\
        &= \sum_{\substack{ |\bar{k} + \bar{l}| \le m \\ \bar{l} \ge \bar{h} }} \frac{ f^{(\bar{k} + \bar{l})}(y) }{ \bar{l}! } \: \frac{\bar{l}!}{(\bar{l} - \bar{h})!} \: (x - y)^{\bar{l} - \bar{h}} \displaybreak[0]\\
        &= \sum_{\substack{ |\bar{k} + \bar{l}| \le m \\ \bar{l} \ge \bar{h} }} \frac{ f^{(\bar{k} + \bar{l})}(y) }{ (\bar{l} - \bar{h})! } \:  \: (x - y)^{\bar{l} - \bar{h}} \\
        &= \sum_{ |\bar{k} + \bar{h} + \bar{l}| \le m } \frac{ f^{(\bar{k} + \bar{h} + \bar{l})}(y) }{ \bar{l}! } \:  \: (x - y)^{\bar{l}} =
        P^{ \bar{k} + \bar{h} }_y (x)
    \end{align*}
    where in the last equality we replaced $\bar{l}$ with $\bar{l} + \bar{h}$.
\end{proof}

\begin{lemma}
    Let $a, b \in F$ and $x \in \RR^n$.
    Then
%    \begin{align}\label{Align_P(x,b)-P(x,a)}
%        P_b(x) - P_a(x) = \sum_{|\bar{l}| \le m} R^{\bar{l}}(b,a) \frac{ (x - b)^{\bar{l}} } { \bar{l}! },
%    \end{align}
%    and more generally,
for each $\bar{k}$ with $|\bar{k}| \le m$,
    \begin{align}\label{Align_P_j(x,b)-P_j(x,a)}
        P^{\bar{k}}_b (x) - P^{\bar{k}}_a (x) = \sum_{|\bar{k} + \bar{l}| \le m} R^{\bar{k} + \bar{l}}(b,a) \frac{ (x - b)^{\bar{l}} } { \bar{l}! }.
    \end{align}
\end{lemma}
\begin{proof}
%    It suffices to show (\ref{Align_P_j(x,b)-P_j(x,a)}) since (\ref{Align_P(x,b)-P(x,a)}) is a special case.
    We expand the polynomial $G_{\bar{k}}(x) := P^{\bar{k}}_b (x) - P^{\bar{k}}_a (x)$ around the point $b \in F$, obtaining
    \begin{align*}
        P^{\bar{k}}_b (x) - P^{\bar{k}}_a (x) &= \sum_{|\bar{k} + \bar{l}| \le m} \frac{ \partial^{\bar{l}} \, G_{\bar{k}} }{ \partial x^{\bar{l}} } \, (b) \: \frac{ (x - b)^{\bar{l}} }{ \bar{l}! }  \\
        &= \sum_{|\bar{k} + \bar{l}| \le m} \left( f^{ ( \bar{k} + \bar{l} ) }(b) -  P^{\bar{k} + \bar{l}}_a (b) \right) \: \frac{ (x - b)^{\bar{l}} }{ \bar{l}! }  \\
        &= \sum_{|\bar{k} + \bar{l}| \le m} R^{\bar{k} + \bar{l}}(b,a)  \: \frac{ (x - b)^{\bar{l}} }{ \bar{l}! },
    \end{align*}
    where in the second equality we used Proposition \ref{Prop_DlPj=Pj+l}, and that for each $\bar{k}$ with $|\bar{k}| \le m$, $P^{\bar{k}}_b (b) = f^{ (\bar{k}) }(b)$.
\end{proof}

We can now define the multi-valued function $\mathrm{WET}_m : \mathcal{J}^m_c(\RR^n) \rightrightarrows C^m(\RR^n)$.
Let $(f^{(\bar{k})})_{|\bar{k}| \le m} \in \mathcal{J}^m_c(\RR^n)$ and let $F$ be its domain.
Given $\FF \in \CUBE (F)$ as in Section \ref{Section_DecompositionCubes} and the partition of unity $\{\varphi^*_Q\}_{Q \in \FF}$ as in Section \ref{Section_PartitionUnit}.
We define $g \in \mathrm{WET}_m ((f^{(\bar{k})})_{|\bar{k}| \le m})$ as a function in $C^m(\RR^n)$ extending $f$ as
\begin{align} \label{Align_Defgj1}
    g(x) =
    \begin{cases}
        f(x) &\textrm{ if }x \in F \\
        \sum_{Q \in \FF} P_{r_Q} (x) \varphi^*_Q(x) &\textrm{ else}.
    \end{cases}
\end{align}
%and, for $|\bar{k}| \le m$,
%\begin{align} \label{Align_Defgj2}
%    g^{(\bar{k})}(x) =
%    \begin{cases}
%        f^{(\bar{k})}(x) &\textrm{ if }x \in F \\
%        \frac{ \partial^{\bar{k}} }{ \partial x^{\bar{k}} } \, g(x) &\textrm{ else}.
%    \end{cases}
%\end{align}
As in Section \ref{Section_FirstWET}, $\mathrm{WET}_m$ is multi-valued since different total names of the closed set $F$ lead to different functions $g$.
Again, if we fix a total name $q \in \baire$ of $F$, we obtain a single valued map $\mathrm{WET}_m^q$ which is linear.

\smallskip
Our first goal in the analysis of $g$ as defined in (\ref{Align_Defgj1}) is to establish upper bounds for the difference between $g^{(\bar k)}$ and the $P_a^{\bar k}$'s.
We first deal with the simplest case, i.e.\ $\bar k = \bar 0$.
%Now, we prove some properties of $g \in \mathrm{WET}_m ((f^{(\bar{k})})_{|\bar{k}| \le m})$.
Let $e := \frac {2}{1 - \varepsilon}$ as in Section \ref{Section_FirstWET}.

\begin{lemma}\label{Lemma_g(x)-P(x)}
    Let $(f^{(\bar{k})})_{|\bar{k}| \le m} \in \mathcal{J}^m(F)$ and suppose that $M$ satisfies the compatibility conditions for the jet.
    Set
    \begin{align}\label{Align_DefA}
        A_m := M \, \sum_{|\bar{l}| \le m} \frac{(7e + 1)^{m - |\bar{l}| + 1}}{ \bar{l}! }.
    \end{align}
    Then, for $g \in \mathrm{WET}_m ((f^{(\bar{k})})_{|\bar{k}| \le m})$ defined as in (\ref{Align_Defgj1}), $x \in \RR^n$ and $a \in F$ we have
    \begin{align} \label{Align_(g(x) - P(x,a))}
        |g(x) - P_a(x)| \le A_m \, d(x,a)^{m + 1}.
    \end{align}
\end{lemma}
\begin{proof}
Notice that $A_m \ge M$.
If $x \in F$, then by Definition \ref{Def_WhitneyJet}
\begin{align*}
    |g(x) - P_a (x)| = |R(x,a)| \le M \, d(x,a)^{m + 1} \le A_m \, d(x,a)^{m + 1}.
\end{align*}

If $x \not\in F$, $g(x) = \sum_{Q \in \FF} P_{r_Q}(x) \varphi^*_Q(x)$. Then, using that $\{ \varphi_Q^* \}_{Q \in \FF}$ is a partition of the unit, (\ref{Align_P_j(x,b)-P_j(x,a)}) for $\bar{k} = \bar{0}$ and (\ref{Align_CompatibilityConditions}) we have
\begin{align*}
    |g(x) - P_a(x)| &= \left\lvert \sum_{Q \in \FF} (P_{r_Q}(x) - P_a(x)) \varphi^*_Q(x) \right\rvert  \\
    &\le \sum_{Q \in \FF_x} \left| \sum_{|\bar{l}| \le m} R^{\bar{l}}(a,r_Q) \frac{ (x - a)^{\bar{l}} } { \bar{l}! } \right| \varphi^*_Q(x)  \\
    &\le \sum_{Q \in \FF_x} \sum_{|\bar{l}| \le m} M \, d(a,r_Q)^{m - |\bar{l}| + 1} \frac{ | (x - a)^{\bar{l}} \, | } { \bar{l}! } \, \varphi^*_Q(x).
\end{align*}
Recall that by Proposition \ref{Prop_(x-y)^l<=d(x,y)^l} $| (x - a)^{\bar{l}} | \le d(x,a)^{|\bar{l}|}$.
Moreover, if $Q \in \FF_x$, we have $d(a,r_Q) \le d(x, a) + d(x, r_Q) \le (7e + 1) \, d(x, a)$ because, by Lemma \ref{Lemma_RelationDistanceApproxProjection}, $d(x,r_Q) \le 7e \, d(x,a)$. %so $d(a, r_Q) \le (7e + 1) \, d(x, a)$, and
Then
\begin{align*}
    |g(x) - P_a(x)| &\le \sum_{Q \in \FF_x} \sum_{|\bar{l}| \le m} M \, (7e + 1)^{m - |\bar{l}| + 1} \, d(x, a)^{m - |\bar{l}| + 1} \: \frac{ d(x, a)^{|\bar{l}|} } { \bar{l}! } \, \varphi^*_Q(x)  \\
     \le \: M \, &d(x,a)^{m + 1} \sum_{Q \in \FF_x} \varphi_Q^*(x) \: \sum_{|\bar{l}| \le m} \frac{(7e + 1)^{m - |\bar{l}| + 1}}{ \bar{l}! } = A_m \, d(x,a)^{m + 1},
\end{align*}
which concludes the proof.
\end{proof}

Now, we want to extend Lemma \ref{Lemma_g(x)-P(x)} to the derivatives of $g$.
To this end, we use Proposition \ref{Prop_Derivativesvarphi}.

\begin{prop}\label{Prop_g^k=sumPk}
    Let $(f^{(\bar{k})})_{|\bar{k}| \le m} \in \mathcal{J}^m(F)$, $g$ be the extension defined in (\ref{Align_Defgj1}) and $x \notin F$.
    For every $\bar{k}$ with $|\bar{k}| \le m$,
    \begin{align*}
        g^{(\bar{k})}(x) = \sum_{Q \in \FF_x} \sum_{\bar{l} \le \bar{k}} \binom{\bar{k}}{\bar{l}} P^{\bar{l}}_{r_Q} (x) \, {\partial_{\bar{k} - \bar{l}}} \, \varphi^*_Q(x) .
    \end{align*}
\end{prop}
\begin{proof}
    % Setting $\delta = d(x, F)/2$\al{$\delta$ non compare mai} and using Proposition \ref{Prop_FydFinito} to interchange limit and sum, we have
    % \begin{align*}
    %     \lim_{h \to 0} \frac{g(x + h \, e_i) - g(x)}{h} &= \lim_{h \to 0} \sum_{Q \in \FF} \frac{ P_{r_Q}(x + h \, e_i) \varphi_Q^*(x + h \, e_i) - P_{r_Q}(x) \varphi_Q^*(x) }{ h } =  \\
    %     &\!\!\!\!\!\!= \lim_{h \to 0} \sum_{Q \in \widetilde{\FF}_x} \frac{ P_{r_Q}(x + h \, e_i) \varphi_Q^*(x + h \, e_i) - P_{r_Q}(x) \varphi_Q^*(x) }{ h } = \displaybreak[0]\\
    %     &\!\!\!\!\!\!=  \sum_{Q \in \widetilde{\FF}_x} \lim_{h \to 0} \frac{ P_{r_Q}(x + h \, e_i) \varphi_Q^*(x + h \, e_i) - P_{r_Q}(x) \varphi_Q^*(x) }{ h } = \\
    %     &\!\!\!\!\!\!=  \sum_{Q \in \widetilde{\FF}_x} \partial _i  ( P_{r_Q} (x) \varphi_Q^*(x) ) = \sum_{Q \in \FF_x} \partial _i  ( P_{r_Q} (x) \varphi_Q^*(x) ).
    % \end{align*}
Proposition \ref{Prop_FydFinito} allows to apply the general Leibniz rule in the ball $B(x,d(x,F)/2)$. Thus, using Proposition \ref{Prop_DlPj=Pj+l}, we obtain
    \begin{align*}
        g^{(\bar{k})}(x) & = \sum_{Q \in \FF_x} \partial_{\bar{k}} \left (P_{r_Q} (x) \varphi^*_Q(x) \right)\\
        & = \sum_{Q \in \widetilde{\FF}_x} \partial_{\bar{k}} \left (P_{r_Q} (x) \varphi^*_Q(x) \right)\\
        & = \sum_{Q \in \widetilde{\FF}_x} \sum_{\bar{l} \le \bar{k}} \binom{\bar{k}}{\bar{l}} P^{\bar{l}}_{r_Q} (x) \, {\partial_{\bar{k} - \bar{l}}} \, \varphi^*_Q(x)\\
        & = \sum_{Q \in \FF_x} \sum_{\bar{l} \le \bar{k}} \binom{\bar{k}}{\bar{l}} P^{\bar{l}}_{r_Q} (x) \, {\partial_{\bar{k} - \bar{l}}} \, \varphi^*_Q(x).\qedhere
    \end{align*}
\end{proof}
 		
Recall from Proposition \ref{Prop_Gx} that for every $x \notin F$, $|\FF_x| \le N_n$.
We let now for every $\bar{k}$ with $|\bar{k}| \le m$
\begin{align}\label{Align_DefA^k}
    A_m^{\bar{k}} &:= \sum_{|\bar{k} + \bar{l}| \le m} M \, \frac{(7e + 1)^{m - |\bar{k} + \bar{l}| + 1}}{ \bar{l} ! }\, + \\ \nonumber
     & \quad\quad + \sum_{\bar{0} < \bar{l} \le \bar{k}} \binom{\bar{k}}{\bar{l}} N_n \, \frac{(98e)^{|\bar{l}|} \, B'_{\bar{l}}}{\varepsilon^{|\bar{l}|}} \sum_{|\bar{k} - \bar{l} + \bar{h} | \le m} M \, \frac{(7e + 1)^{m - |\bar{k}-\bar{l}+\bar{h}| + 1}}{\bar{h}!}.
\end{align}
Notice that $A_m^{\bar{0}} = A_m$ defined in (\ref{Align_DefA}).
 		
\begin{lemma}\label{Lemma_g^k-P^k}
    Let $(f^{(\bar{k})})_{|\bar{k}| \le m} \in \mathcal{J}^m(F)$ and $g$ be the extension defined in (\ref{Align_Defgj1}).
    Then, for all $x \in \RR^n$, $a \in F$ with $d(x,a) \le 7e \, d(x, F)$, and for each $\bar{k}$ with $|\bar{k}| 	\le m$,
    \begin{align} \label{Align_(g^k(x) - P^k(x,a))}
        |g^{(\bar{k})}(x) - P^{\bar{k}}_a(x)| \le A_m^{\bar{k}} \, d(x,a)^{m - |\bar{k}| + 1}.
    \end{align}
\end{lemma}
\begin{proof}
    If $x \in F$, then by Definition \ref{Def_WhitneyJet}
    \begin{align*}
        |g^{(\bar{k})}(x) - P^{\bar{k}}_a(x)| \le M \, d(x,a)^{m - |\bar{k}| + 1} \le A_m^{\bar{k}} \, d(x,a)^{m - |\bar{k}| + 1}.
    \end{align*}

    If instead $x \not\in F$, using Propositions \ref{Prop_g^k=sumPk}, \ref{Prop_DlPj=Pj+l} and \ref{Prop_varphi*Partition} we have
    \begin{multline*}
        |g^{(\bar{k})}(x) - P^{\bar{k}}_a(x)| = \left| \sum_{Q \in \FF_x} \sum_{\bar{l} \le \bar{k}} \binom{\bar{k}}{\bar{l}} P^{\bar{l}}_{r_Q} (x) \, {\partial_{\bar{k} - \bar{l}}} \, \varphi^*_Q(x) - P^{\bar{k}}_a(x) \right| \\
        = \left| \sum_{Q \in \FF_x} \left( P^{\bar{k}}_{r_Q} (x) - P^{\bar{k}}_a(x) \right) \varphi^*_Q(x) + \sum_{\bar{0} < \bar{l} \le \bar{k}} \binom{\bar{k}}{\bar{l}} \sum_{Q \in \FF_x}  P^{\bar{k}-\bar{l}}_{r_Q}(x) \, {\partial_{\bar{l}}} \varphi^*_Q(x) \right|.
    \end{multline*}
    For the first term of the sum, arguing as in the proof of Lemma \ref{Lemma_g(x)-P(x)}, we have
    \begin{align}\label{Align_LemmaAkm}
        \nonumber &\!\!\!\!\!\!\!\!\!\!\!\!\!\!\!\!\!\!\left| \sum_{Q \in \FF_x} ( P^{\bar{k}}_{r_Q} (x) - P^{\bar{k}}_a(x) ) \varphi^*_Q(x) \right| \le \\
        &\le \sum_{Q \in \FF_x}  \sum_{|\bar{k} + \bar{l}| \le m} M \, d(a, r_Q)^{m - |\bar{k} + \bar{l}| + 1} \frac{ d(x,a)^{|\bar{l}|} }{ \bar{l} ! } \varphi^*_Q(x)  \\
        \nonumber &\le \sum_{Q \in \FF_x}  \sum_{|\bar{k} + \bar{l}| \le m} M \, \frac{(7e + 1)^{m - |\bar{k} + \bar{l}| + 1}}{ \bar{l} ! } \, d(x, a)^{m - |\bar{k}| + 1} \, \varphi^*_Q(x).
    \end{align}

    Now, let $\bar{l}$ with $\bar{0} < \bar{l} \le \bar{k}$ be fixed. By Proposition \ref{Prop_varphi*Partition} it follows that for every $x \notin F$,
    \begin{align}\label{Align_sumvarphider0}
         \sum_{Q \in \FF_x} {\partial_{\bar{l}}} \varphi^*_Q(x) = 0.
    \end{align}
    Therefore,
    \begin{align*}
        \sum_{Q \in \FF_x}  P^{\bar{k}-\bar{l}}_{r_Q}(x) \, {\partial_{\bar{l}}} \varphi^*_Q(x) = \sum_{Q \in \FF_x}  \left( P^{\bar{k}-\bar{l}}_{r_Q}(x) - P^{\bar{k}-\bar{l}}_{a}(x) \right) {\partial_{\bar{l}}} \varphi^*_Q(x).
    \end{align*}
    Using (\ref{Align_NormVarphi^*DiamQ}), $d(x,a) \le 7e \,d(x,F)$, and $\diam(Q) > d(x,F) / 7$ (as shown in the proof of Proposition \ref{Prop_Gx}) we have
    \begin{align} \label{Align_|dephiQ|le49e}
        \left| {\partial_{\bar{l}}} \varphi^*_Q(x) \right| \le B'_{\bar{l}} \left( \frac2{\varepsilon \, \diam(Q) }\right) ^{|\bar{l}|}  \le  B'_{\bar{l}} \left( \frac{ 98e } {\varepsilon \, d(x,a) } \right) ^{|\bar{l}|}.
    \end{align}
    Adapting (\ref{Align_LemmaAkm}), we then have
    \begin{align*}
        &\left| \sum_{Q \in \FF_x}  P^{\bar{k}-\bar{l}}_{r_Q}(x) \, {\partial_{\bar{l}}} \varphi^*_Q(x) \right| = \left| \sum_{Q \in \FF_x}  \left( P^{\bar{k}-\bar{l}}_{r_Q}(x) - P^{\bar{k}-\bar{l}}_{a}(x) \right) {\partial_{\bar{l}}} \varphi^*_Q(x) \right|  \\
        &\le \sum_{Q \in \FF_x}  \left| P^{\bar{k}-\bar{l}}_{r_Q}(x) - P^{\bar{k}-\bar{l}}_{a}(x) \right|  \left| {\partial_{\bar{l}}} \varphi^*_Q(x) \right| \displaybreak[0]\\
        &\le \sum_{Q \in \FF_x} \left( \sum_{|\bar{k} - \bar{l} + \bar{h} | \le m} M \, \frac{(7e + 1)^{m - |\bar{k}-\bar{l}+\bar{h}| + 1}}{\bar{h}!} \, d(x,a)^{m - |\bar{k}-\bar{l}| + 1} \right) \, B'_{\bar{l}} \left( \frac{  98e } {\varepsilon \, d(x,a) } \right) ^{|\bar{l}|} \\
        &\le \left( N_n \, \frac{ (98e)^{|\bar{l}|} \, B'_{\bar{l}} }{\varepsilon^{|\bar{l}|} } \sum_{|\bar{k} - \bar{l} + \bar{h} | \le m} M \, \frac{(7e + 1)^{m - |\bar{k}-\bar{l}+\bar{h}| + 1}}{\bar{h}!} \, \right) \,  d(x,a)^{m-|\bar{k}|+1}.
    \end{align*}
    Now, putting all together and using again Proposition \ref{Prop_varphi*Partition}, we obtain
    \begin{align*}
        &|g^{(\bar{k})}(x) - P^{\bar{k}}_a(x)| \le \\
        &\left| \sum_{Q \in \FF_x} \left( P^{\bar{k}}_{r_Q} (x) - P^{\bar{k}}_a(x) \right) \varphi^*_Q(x) \right| + \sum_{\bar{0} < \bar{l} \le \bar{k}} \binom{\bar{k}}{\bar{l}} \left| \sum_{Q \in \FF_x}  P^{\bar{k}-\bar{l}}_{r_Q}(x) \,{\partial_{\bar{l}}} \varphi^*_Q(x) \right|  \displaybreak[0]\\
        &\le \sum_{|\bar{k} + \bar{l}| \le m} M \, \frac{(7e + 1)^{m - |\bar{k} + \bar{l}| + 1}}{ \bar{l} ! } \, d(x,a)^{m-|\bar{k}|+1} \,+ \\
        &\quad\quad \sum_{\bar{0} < \bar{l} \le \bar{k}} \binom{\bar{k}}{\bar{l}} N_n \, \frac{ (98e)^{|\bar{l}|}  \, B'_{\bar{l}}}{\varepsilon^{|\bar{l}|}} \sum_{|\bar{k} - \bar{l} + \bar{h} | \le m} M \, \frac{(7e + 1)^{m - |\bar{k}-\bar{l}+\bar{h}| + 1}}{\bar{h}!} \, d(x,a)^{m-|\bar{k}|+1} \\
        &= A^{\bar{k}}_m \, d(x,a)^{m-|\bar{k}|+1}.\qedhere
    \end{align*}
    %Therefore the desired result is achieved.
\end{proof}

\begin{prop}
    Let $F \subseteq \RR^n$ be closed, and let $(f^{(\bar{k})})_{|\bar{k}| \le m} \in \mathcal{J}^m_c(\RR^n)$ be a Whitney jet of order $m$ on $F$.
    For every $g \in \mathrm{WET}_m((f^{(\bar{k})})_{|\bar{k}| \le m})$ the following hold:
    \begin{enumerate}
%        \item\label{Item_WET_1} $g(x) = f^{(\bar{0})}(x)$ on $F$,
        \item\label{Item_WET_2} ${\partial_{\bar{k}}} \,g(x) = f^{(\bar{k})}(x)$ on $F$ for all $\bar{k}$ with $|\bar{k}| \le m$,
        \item\label{Item_WET_3} $g$ is $C^\infty$ on $\RR^n \setminus F$.
    \end{enumerate}
\end{prop}
\begin{proof}
    We prove \ref{Item_WET_2} and \ref{Item_WET_3} by induction on $m$.
    For $m = 0$, condition \ref{Item_WET_2} follows from the properties of $\mathrm{WET}_0$ given in the previous section.
    For \ref{Item_WET_3}, notice that on every point of $\RR^n \setminus F$, $g$ is a finite sum of $C^{\infty}$ functions (see Remark \ref{Cor_SumVarphi} and Proposition \ref{Prop_varphiCinfinito}).

    For the case $m = 1$, we have to show that $g$ is a continuous extension of $f$, and also that $g$ has continuous partial derivatives on $\RR^n$ which coincide with $f^{(\bar{k})}$ on $F$.

    First, we prove that $g$ is continuous by showing that for all $y \in \RR^n$, $\lim_{x \to y} g(x) = g(y)$.

    If $y \in F$, we use inequality (\ref{Align_(g(x) - P(x,a))}) with $y = a$, and obtain
    $$ \lim_{x \to y} |g(x) - P_y(x)| \le \lim_{x \to y} A_1 \, d(x,y)^2 = 0.$$
    So $\lim_{x \to y} g(x) = \lim_{x \to y} P_y(x) = P_y(y) = f(y) = g(y)$ by (\ref{Align_P(x,y)}), and (\ref{Align_Defgj1}).

     If $y \notin F$, setting $\delta = d(y,F)/2$ and using Proposition \ref{Prop_FydFinito} to interchange limit and sum in $B(y,\delta)$, we have
    \begin{align*}
        \left| g(y) - \lim_{x \to y} g(x) \right| &= \lim_{x \to y} \left | \sum_{Q \in \FF} P_{r_Q} (y) \varphi_Q^*(y) -  \sum_{Q \in \FF} P_{r_Q} (x) \varphi_Q^*(x) \right|  \\
        &\le \lim_{x \to y} \sum_{Q \in \FF} \left| P_{r_Q} (y) \varphi_Q^*(y) - P_{r_Q} (x) \varphi_Q^*(x) \right| \displaybreak[0]\\
        &= \lim_{x \to y} \sum_{Q \in \widetilde{\FF}_y} \left| P_{r_Q} (y) \varphi_Q^*(y) - P_{r_Q} (x) \varphi_Q^*(x) \right| \\
        &=  \sum_{Q \in \widetilde{\FF}_y} \lim_{x \to y} | P_{r_Q} (y) \varphi_Q^*(y) - P_{r_Q} (x) \varphi_Q^*(x) | = 0,
    \end{align*}
    since both $P_{r_Q}$ and $\varphi_Q^*$ are continuous, and therefore also their product.

    Next, we prove that $g$ has partial derivatives at every point of $\RR^n$, that is, that for every $y \in \RR^n$ and for all $1 \le i \le n$, the limit
    $$ \lim_{h \to 0} \frac{ g(y + h \, e_i) - g(y) }{h},$$
    exists and is finite, and these partial derivatives coincide with $f^{(\bar{k})}$ on $F$, and are continuous in $\RR^n$.

    Let $y \in F$, and for $1 \le i \le n$, we call $\bar{\imath}$ the vector $e_i$ seen as a $n$-tuple in $\NN$.
    We have
    \begin{align*}
        \left| \lim_{h \to 0} \frac{ g(y + h \, e_i) - g(y) }{h} - f^{(\bar{\imath})}(y) \right| &= \lim_{h \to 0} \left| \frac{ \, g(y + h \, e_i) - f(y) - h \, f^{(\bar{\imath})}(y)}{h} \right|  \\
        &= \lim_{h \to 0} \left| \frac{ g(y + h \, e_i) - P_y( y + h \, e_i)} {h} \right|  \\
        & \le  \lim_{h \to 0} \frac{1}{|h|} \, A_1 \, d(y + h \, e_i, y)^2 = \lim_{h \to 0} A_1 \, |h| = 0.
    \end{align*}

    \noindent Therefore, the partial derivatives of $g$ in $y$ exist and coincide with $f^{(\bar{\imath})}(y)$, that is
    \begin{align}
        \partial_{\bar{\imath}} g(y) = f^{(\bar{\imath})}(y) \quad \text{for all }y \in F.
    \end{align}

    \noindent Next, we prove that they are continuous on $\RR^n$. Since $m = 1$, given $\bar{k}$ with $|\bar{k}|=1$, we have $P^{\bar{k}}_y (x) = f^{(\bar{k})}(y)$.
    Then, using inequality (\ref{Align_(g^k(x) - P^k(x,a))}) we obtain
    $$ \lim_{x \to y} |g^{(\bar{k})} (x) - f^{(\bar{k})}(y)| = \lim_{x \to y} |g^{(\bar{k})}(x) - P^{\bar{k}}_y (x)| \le \lim_{x \to y} A_m^{\bar{k}} \, d(x,y) = 0,$$
    so we have
    \begin{align} \label{Align_g^kContinuousF}		
        \lim_{x \to y} g^{(\bar{k})}(x) = g^{(\bar{k})}(y) \quad \text{ for all }y \in F.
    \end{align}
   If instead $y \notin F$, using Proposition \ref{Prop_g^k=sumPk}, and that $P_{r_Q}$ and $\varphi_Q^*$ are $C^{\infty}$ (Proposition \ref{Prop_varphiCinfinito}), it follows that $g$ is differentiable at $y$.

    It remains to prove that the partial derivatives of $g$ are continuous on $\RR^n$, but this simply follows from (\ref{Align_g^kContinuousF}) and Proposition \ref{Prop_g^k=sumPk}.
    Therefore, we have proved that the extension of $(f^{(\bar{k})})$ defined in (\ref{Align_Defgj1}) satisfies condition \ref{Item_WET_2}. %conditions \ref{Item_WET_1} and \ref{Item_WET_2}.

    For \ref{Item_WET_3}, it suffices to observe that, by Proposition \ref{Prop_g^k=sumPk}, not only $g$ is continuously differentiable on $\RR^n \setminus F$, but also $C^\infty$.

    The proof for $m \ge 2$ can then be carried out by induction, the inductive step being very similar to the case $m = 1$ just given.
\end{proof}

\begin{teo}[Computable Whitney Extension Theorem]\label{Theorem_WETm}
    For all $m \ge 0$ the multi-valued function $\mathrm{WET}_m$ is computable.
\end{teo}
\begin{proof}
We show that for each $q$ which is a $\psi$-name for the closed set $F$ the map $\mathrm{WET}_m^q$ is computable.
We fix once for all $q$ and hence also $\FF$, $\{ \varphi_Q^* \mid Q \in \FF \}$ and $r_Q$ for all $Q \in \FF$ and argue by induction on $m$.

The case $m = 0$ is the First Computable Whitney Extension Theorem (Theorem \ref{Theorem_FirstWET}).

Now assume that $\mathrm{WET}^q_k$ is computable for all $k \le m$ and we are given a Whitney jet $J_{m+1} := (f^{(\bar{l})})_{|\bar{l}| \le m+1}$.
Let $J_m$ be the restriction of $J_{m+1}$ up to order $m$, and $g_m = \mathrm{WET}_m^q(J_m)$ (recall that $g_m$ is actually encoded as $(g_m^{(\bar{l})})$). We now show how to compute $g_{m+1} = \mathrm{WET}^q_{m+1}(J_{m+1})$, encoded as $(g_{m+1}^{(\bar{l})})$.
For every $y \in F$, we call $P^{\bar{k}}_{y,m}(x)$ the polynomials defined in (\ref{Align_P(x,y)}) corresponding to $J_m$, and $P^{\bar{k}}_{y,m+1}(x)$ the analogous polynomials corresponding to $J_{m+1}$.
It follows from the definition that for every $\bar{k}$
    \begin{align}\label{Align_P_ym_P_ym+1}
        P_{y, m+1}^{\bar{k}}(x) = P_{y, m}^{\bar{k}}(x) + \sum_{|\bar{k} + \bar{l}| = m+1} \frac{1}{\bar{l}!} f^{(\bar{k} + \bar{l})}(y) \, (x - y)^{\bar{l}}.
    \end{align}
    Let also $p(m) := \left |\left\{ \bar{l} \in \NN^n \mid |\bar{l}| = m \right\} \right|$.

    Now, as usual, we are given $x \in \RR^n$.
    Let us start from the computation of $g_{m+1}$ as $g^{(\bar{0})}_{m+1}$.
    We apply $\mathrm{WET}_0^q$ to the functions $f^{(\bar{l})}$, with $|\bar{l}| = m+1$, and compute an approximation $s_{\bar{l}}$ of $h^{\bar{l}}_0(x) =: \mathrm{WET}_0^q(f^{(\bar{l})})(x)$ up to $\frac{1}{2}$.
    Then we set
    $$ S := \frac{1}{2} + \max_{|\bar{l}| = m+1} \frac{1}{\bar{l}!} \, |s_{\bar{l}}|,$$
    and $c := 14 e + 1$ (as in the proof of Theorem \ref{Theorem_FirstWET}).

    At stage $i$ we search simultaneously for an open basic ball $B(y, \delta)$ in the negative information of $F$ with $x \in B(y, \delta)$ (in this case we are sure that $x \not\in F$ and we can compute $g_{m+1}(x)$ using the definition given in (\ref{Align_Defgj1})) and for a pair of the form $(B(y'_{\bar{0}}, \delta), B(y''_{\bar{0}}, \epsilon_{\bar{0}}))$ in the information of $g_m$ for such that
    \begin{enumerate}
        \item\label{Item_g_m+1_1} $B(y'_{\bar{0}}, \delta)$ is in the positive information of $F$ and $x \in B(y'_{\bar{0}}, \delta)$,
        \item $\epsilon_{\bar{0}} < 2^{-i-1}$,
        \item $p(m+1) \, S \, ( c \, \delta )^{m+1} < 2^{-i-1}$.
    \end{enumerate}
    If we find such a pair, we claim that it is correct to write $y''_{\bar{0}}$ on the output tape.

    In fact, if $x \in F$, then
    \begin{align*}
        |y''_{\bar{0}} - f(x)| < \epsilon_{\bar{0}} < 2^{-i}.
    \end{align*}
    If instead $x \notin F$ let $y \in B(y'_{\bar{0}}, \delta) \cap F$ (such a $y$ exists by condition (\ref{Item_g_m+1_1})) and notice that when $Q \in \FF_x$, by Lemma \ref{Lemma_RelationDistanceApproxProjection} we have $d(x,r_Q) < 7e \, d(x,y) < 14e \, \delta < c \, \delta$.
    Using (\ref{Align_P_ym_P_ym+1}) with $\bar{k} = \bar{0}$ and Proposition \ref{Prop_(x-y)^l<=d(x,y)^l} we obtain
    \begin{align*}
        |y''_{\bar{0}} - g_{m+1}(x)| &\le \left| y''_{\bar{0}} - \sum_{Q \in \FF_x} P_{r_Q, m+1}(x) \varphi_Q^*(x) \right|  \\
        &\le | y''_{\bar{0}} - g_{m}(x) | +  \sum_{Q \in \FF_x} \left| \sum_{|\bar{l}| = m+1} \frac{1}{\bar{l}!} f^{(\bar{l})}(r_Q) \, (x - r_Q)^{\bar{l}} \right| \varphi_Q^*(x)  \\
        &< \epsilon_{\bar{0}} + \sum_{Q \in \FF_x} \left( \, p(m+1) \, S \, ( c \, \delta )^{m+1} \right)\, \varphi_Q^*(x) < 2^{-i}
    \end{align*}
    as $\frac{1}{\bar{l}!} f^{(\bar{l})}(r_Q)\leq S$ is guaranteed by the proof of Theorem \ref{Theorem_FirstWET} for the case $i = 1$ where $y'' = h_0^{\bar{l}}(x)[1]$, because $|y'' - f^{(\bar{l})}(r_Q)| < 1/2$.

    If the pair $(B(y'_{\bar{0}}, \delta), B(y''_{\bar{0}}, \epsilon_{\bar{0}}))$ does not exist, then we will eventually realize that $x \notin F$, and when this happens we compute an approximation of $g_m(x)$ as explained above.\smallskip

    To compute $g^{(\bar{l})}_{m+1}(x)$ with $|\bar{l}| >0$, we start with the case $|\bar{l}| \le m$.
    Using the definition in (\ref{Align_P(x,y)}) and Proposition \ref{Prop_g^k=sumPk}, when $x \notin F$, we have
    \begin{align}
        g^{(\bar{l})}_{m}(x) &= \sum_{Q \in \FF_x} \sum_{ \bar{h} \le \bar{l}} \binom{\bar{l}}{\bar{h}} P_{r_Q, m}^{\bar{h}}(x) \, {\partial _{\bar{l}-\bar{h}}} \, \varphi_Q^*(x), \\
        \label{Align_g_m+1l} g^{(\bar{l})}_{m+1}(x) &= \sum_{Q \in \FF_x} \sum_{\bar{h} \le \bar{l}} \binom{\bar{l}}{\bar{h}} P_{r_Q, m+1}^{\bar{h}}(x) \, {\partial _{\bar{l}-\bar{h}}} \,  \varphi_Q^*(x).
    \end{align}
    Using (\ref{Align_P_ym_P_ym+1}) and Propositions \ref{Prop_DlPj=Pj+l}, \ref{Prop_g^k=sumPk}, we obtain
    \begin{align*}
        g^{(\bar{l})}_{m+1}(x) &= \sum_{Q \in \FF_x} \sum_{ \bar{h} \le \bar{l}} \binom{\bar{l}}{\bar{h}} P_{r_Q, m}^{\bar{h}}(x) \, {\partial _{\bar{l}-\bar{h}}} \,  \varphi_Q^*(x) + \\
        & \sum_{Q \in \FF_x} \left( \sum_{ \bar{h} \le \bar{l}} \binom{\bar{l}}{\bar{h}} \, {\partial _{\bar{l}-\bar{h}}} \,  \varphi_Q^*(x) \sum_{|\bar{h} + \bar{k}| = m+1}  \frac{1}{\bar{k}!} \, f^{(\bar{h} + \bar{k})}(r_Q) \, (x - r_Q)^{\bar{k}} \right) \\
        &\!\!\!\!\!\!\!\!\!\!\!\!\!\!\!\! = g^{(\bar{l})}_{m}(x) + \sum_{Q \in \FF_x} \left( \sum_{\bar{h} \le \bar{l}} \binom{\bar{l}}{\bar{h}} \, {\partial _{\bar{l}-\bar{h}}} \,  \varphi_Q^*(x) \sum_{|\bar{h} + \bar{k}| = m+1}  \frac{1}{\bar{k}!} \, f^{(\bar{h} + \bar{k})}(r_Q) \, (x - r_Q)^{\bar{k}} \right).
    \end{align*}
    Now, as for the case $\bar{l} = \bar{0}$, we have all the approximations $s_{\bar{r}}$ of $\mathrm{WET}_0^q(f^{(\bar{r})})(x)$ for $|\bar{r}| = m+1$.
    \begin{align*}
        S_{\bar{h}} := \frac{1}{2} + \max_{|\bar{h} + \bar{k}| = m+1} \frac{1}{\bar{k}!} \, |s_{|\bar{h} + \bar{k}|}|,
    \end{align*}
    and define $H_{\bar{l}}$ as
    \begin{align*}
        H_{\bar{l}} := N_n \, \sum_{\bar{h} \le \bar{l}} \binom{\bar{l}}{\bar{h}} \frac{ S_{\bar{h}}\, (98e)^{|\bar{l} - \bar{h}|} \, B'_{\bar{l} - \bar{h}} \, p_{\bar{h} }}{\varepsilon^{|\bar{l} - \bar{h}|}}.
    \end{align*}

    At stage $i$ we search simultaneously for information ensuring that $x \not\in F$ (in which case we compute $g_{m+1}^{(\bar{l})}(x)$ using (\ref{Align_g_m+1l})) and for a pair of the form $(B(y'_{\bar{l}}, \delta), B(y''_{\bar{l}}, \epsilon_{\bar{l}}))$ in the information about $g_{m}^{(\bar{l})}$ such that
    \begin{enumerate}
        \item $B(y'_{\bar{l}}, \delta)$ is in the positive information of $F$ and $x \in B(y'_{\bar{l}}, \delta)$,
        \item $\epsilon_{\bar{l}} < 2^{-i-1}$,
        \item $H_{\bar{l}} \, (c \, \delta)^{m+1 - |\bar{l}|} < 2^{-i-1}$.
    \end{enumerate}
    If we find a pair that meets the requirements, we claim that it is correct to write $y''_{\bar{l}}$ on the output tape.
    In fact, if $x \in F$, using that $g_{m}^{(\bar{l})}(\overline{B(y'_{\bar{l}}, \delta)}) \subseteq B(y''_{\bar{l}}, \epsilon_{\bar{l}}))$, we have
    $$|y''_{\bar{l}} - f^{(\bar{l})}(x)| = |y''_{\bar{l}} - g_m^{(\bar{l})}(x)| \le \epsilon_{\bar{l}} < 2^{-i}.$$

    If $x \not\in F$, we use
    \begin{align*}
        |y''_{\bar{l}} - g_{m+1}^{(\bar{l})}(x)| \le |y''_{\bar{l}} - g_{m}^{(\bar{l})}(x)| + |g_{m+1}^{(\bar{l})}(x) - g_{m}^{(\bar{l})}(x)|.
    \end{align*}
    For the first term on the right, as for the case $x \in F$, we have $|y''_{\bar{l}} - g_{m}^{(\bar{l})}(x)| \le \epsilon_{\bar{l}}  < 2^{-i-1}$.

    For the second term, notice that, using a similar argument to that in the proof of Lemma \ref{Lemma_g^k-P^k}, we obtain, for every $Q \in \FF_x$,
   \begin{align}\label{Align_varphi*Qle49e}
         \left| {\partial _{\bar{l}-\bar{h}}} \, \varphi_{Q}^*(x) \right| \le \frac{ {(98e)^{|\bar{l} - \bar{h}|}} \, B'_{\bar{l}-\bar{h}} }{ \varepsilon^{|\bar{l}-\bar{h}|} \, d(x, r_Q) ^{|\bar{l}-\bar{h}|}}.
    \end{align}
    Using (\ref{Align_Defgj1}), (\ref{Align_P_ym_P_ym+1}), (\ref{Align_varphi*Qle49e}) and Proposition \ref{Prop_(x-y)^l<=d(x,y)^l}, we then have
    \begin{align*}
        |g_{m+1}^{(\bar{l})}(x) & - g_{m}^{(\bar{l})}(x)| = \\
        &\left| \sum_{Q \in \FF_x} \sum_{\bar{h} \le \bar{l}} \binom{\bar{l}}{\bar{h}} \, {\partial _{\bar{l}-\bar{h}}} \,  \varphi_Q^*(x) \sum_{|\bar{h} + \bar{k}| = m+1}  \frac{1}{\bar{k}!} \, f^{(\bar{h} + \bar{k})}(r_Q) \, (x - r_Q)^{\bar{k}} \right|  \displaybreak[0]\\
        &\le \sum_{Q \in \FF_x} \sum_{ \bar{h} \le \bar{l}} \binom{\bar{l}}{\bar{h}} \left| \, {\partial _{\bar{l}-\bar{h}}} \, \varphi_Q^*(x) \right| \sum_{|\bar{h} + \bar{k}| = m+1}  \left| \frac{1}{\bar{k}!} \, f^{(\bar{h} + \bar{k})}(r_Q) \right| \, d(x,r_Q)^{|\bar{k}|} \displaybreak[0]\\
        &\le \sum_{Q \in \FF_x} \sum_{ \bar{h} \le \bar{l}} \binom{\bar{l}}{\bar{h}} S_{\bar{h}} \frac{ { (98e)^{|\bar{l} - \bar{h}|} } \, B'_{\bar{l} - \bar{h}} }{\varepsilon^{|\bar{l} - \bar{h}|}} \, d(x,r_Q)^{|\bar{h}|-|\bar{l}|} \, p_{\bar{h}} \, d(x,r_Q)^{m+1 - |\bar{h}|}  \displaybreak[0]\\
        &\le \sum_{Q \in \FF_x} \sum_{ \bar{h} \le \bar{l}} \binom{\bar{l}}{\bar{h}} \frac{ S_{\bar{h}}\, { (98e)^{|\bar{l} - \bar{h}|} } \, B'_{\bar{l} - \bar{h}} \, p_{\bar{h} }}{\varepsilon^{|\bar{l} - \bar{h}|}} \, d(x,r_Q)^{m+1 - |\bar{l}|}  \\
        &\le N_n (c \, \delta)^{m+1 - |\bar{l}|} \, \sum_{ \bar{h} \le \bar{l}} \binom{\bar{l}}{\bar{h}} \frac{ S_{\bar{h}}\, { (98e)^{|\bar{l} - \bar{h}|}  } \, B'_{\bar{l} - \bar{h}} \, p_{\bar{h} }}{\varepsilon^{|\bar{l} - \bar{h}|}} = H_{\bar{l}} (c \, \delta)^{m+1 - |\bar{l}|}.
    \end{align*}
    Therefore
    \begin{align*}
        |g_{m+1}^{(\bar{l})}(x) - g_{m}^{(\bar{l})}(x)| \le H_{\bar{l}} (c \, \delta)^{m+1 - |\bar{l}|} < 2^{-i-1},
    \end{align*}
    and $y''_{\bar{l}}$ is a good approximation of $g_{m+1}^{(\bar{l})}(x)$.

    If the pair $(B(y'_{\bar{l}}, \delta), B(y''_{\bar{l}}, \epsilon_{\bar{l}}))$ does not exist, then we will eventually realize that $x \notin F$, and when this happens we compute an approximation of $g_m(x)$ as explained above.\smallskip

    It remains to show how to compute $g_{m+1}^{(\bar{l})}(x)$ when $|\bar{l}| = m+1$. To do this, we cannot use $g_{m}^{(\bar{l})}$ because it does not exist.
    Let $h_0^{\bar{l}} = \mathrm{WET}^q_0(f^{(\bar{l})})$; by (\ref{Align_g_m+1l}), for $x\notin F$ we have
    \begin{align}
    \begin{split}\label{Align_g_m+1}
        g^{(\bar{l})}_{m+1}(x) &= \sum_{Q \in \FF_x} \sum_{ \bar{h} \le \bar{l}} \binom{\bar{l}}{\bar{h}} P_{r_Q, m+1}^{\bar{h}}(x) \, {\partial _{\bar{l}-\bar{h}}} \, \varphi_Q^*(x)  \\
        &= \sum_{Q \in \FF_x} f^{(\bar{l})}(r_Q) \, \varphi_Q^*(x) + \sum_{ \bar{h} < \bar{l}} \binom{\bar{l}}{\bar{h}} \sum_{Q \in \FF_x} P_{r_Q, m+1}^{\bar{h}}(x) \, {\partial _{\bar{l}-\bar{h}}} \, \varphi_Q^*(x)  \\
        &= h_0^{\bar{l}}(x) + \sum_{\bar{h} < \bar{l}} \binom{\bar{l}}{\bar{h}} \sum_{Q \in \FF_x} P_{r_Q, m+1}^{\bar{h}}(x) \, {\partial _{\bar{l}-\bar{h}}} \, \varphi_Q^*(x).
    \end{split}
    \end{align}
    We first define $H_{\bar{l}}$ to be the number
    \begin{align*}
        H_{\bar{l}} := N_n \, \sum_{ \bar{h} < \bar{l}} \binom{\bar{l}}{\bar{h}} \frac{ { (98e)^{|\bar{l} - \bar{h}|} } \, B'_{\bar{l}-\bar{h}} \, A_{m+1}^{\bar{h}} }{ \varepsilon^{|\bar{l}-\bar{h}|}}.
    \end{align*}

    At stage $i$, as usual, we search simultaneously for information ensuring that $x \not\in F$ (in which case we compute $g_{m+1}^{(\bar{l})}(x)$ using (\ref{Align_g_m+1l})) and for a pair of the form $(B(y'_{\bar{l}}, \delta), B(y''_{\bar{l}}, \epsilon_{\bar{l}}))$ in the information about $h_{0}^{\bar{l}}$ such that
    \begin{enumerate}
        \item $B(y'_{\bar{l}}, \delta)$ is in the positive information of $F$ and $x \in B(y'_{\bar{l}}, \delta)$,
        \item $\epsilon_{\bar{l}} < 2^{-i-1}$,
        \item $H_{\bar{l}} \, c \, \delta < 2^{-i-1}$.
    \end{enumerate}
If we find a pair that meets the requirements, we print $y''_{\bar{l}}$ on the output tape.
    In fact, if $x \in F$, we want to have $g^{(\bar{l})}_{m+1}(x)=f^{(\bar{l})}(x)=h_0^{(\bar{l})}(x)$, and hence $|y''_{\bar{l}} - g_{m+1}^{(\bar{l})}(x)| < 2^{-i}$.

    If $x \not\in F$, using formula (\ref{Align_g_m+1}), we have
    \begin{align*}
        |y''_{\bar{l}} - g_{m+1}^{(\bar{l})}(x)| \le |y''_{\bar{l}} - h_{0}^{\bar{l}}(x)| + \left| \sum_{ \bar{h} < \bar{l}} \binom{\bar{l}}{\bar{h}} \sum_{Q \in \FF_x} P_{r_Q, m+1}^{\bar{h}}(x) \,{\partial _{\bar{l}-\bar{h}}}\, \varphi_Q^*(x) \right|.
    \end{align*}
    The first term on the right is less than $2^{-i-1}$. For the second term, fix $\bar{h}$ with $ \bar{h} < \bar{l}$. Using (\ref{Align_sumvarphider0}), Lemma \ref{Lemma_g^k-P^k} and (\ref{Align_varphi*Qle49e}) we have
    \begin{align*}
        \left| \sum_{Q \in \FF_x} P_{r_Q, m+1}^{\bar{h}}(x) \,{\partial _{\bar{l}-\bar{h}}}\, \varphi_Q^*(x) \right| & = \left| \sum_{Q \in \FF_x} \left( P_{r_Q, m+1}^{\bar{h}}(x) - g_{m+1}^{(\bar{h})}(x) \right) \,{\partial _{\bar{l}-\bar{h}}}\, \varphi_Q^*(x) \right| \\
        &\le \sum_{Q \in \FF_x} \left| P_{r_Q, m+1}^{\bar{h}}(x) - g_{m+1}^{(\bar{h})}(x) \right| \, \left| {\partial _{\bar{l}-\bar{h}}}\, \varphi_Q^*(x) \right|  \displaybreak[0]\\
        &\le \sum_{Q \in \FF_x} A_{m+1}^{\bar{h}} \, d(x, r_Q)^{m+2-|\bar{h}|} \, \frac{ { (98e)^{|\bar{l} - \bar{h}|} } \, B'_{\bar{l}-\bar{h}} }{ \varepsilon^{|\bar{l}-\bar{h}|} \, d(x, r_Q) ^{|\bar{l}|-|\bar{h}|}}  \\
        &\le \frac{ { (98e)^{|\bar{l} - \bar{h}|} } \, B'_{\bar{l}-\bar{h}} \, A_{m+1}^{\bar{h}} }{ \varepsilon^{|\bar{l}-\bar{h}|}} \, \sum_{Q \in \FF_x}  \, d(x, r_Q).
    \end{align*}
    Therefore, by the third condition, we have
    \begin{align*}
        &\!\!\!\!\!\!\!\!\!\left| \sum_{ \bar{h} < \bar{l}} \binom{\bar{l}}{\bar{h}} \sum_{Q \in \FF_x} P_{r_Q, m+1}^{\bar{h}}(x) \, {\partial _{\bar{l}-\bar{h}}} \, \varphi_Q^*(x) \right|  \\
        &\le \sum_{ \bar{h} < \bar{l}} \binom{\bar{l}}{\bar{h}} \frac{ { (98e)^{|\bar{l} - \bar{h}|} } \, B'_{\bar{l}-\bar{h}} \, A_{m+1}^{\bar{h}} }{ \varepsilon^{|\bar{l}-\bar{h}|}} \, \sum_{Q \in \FF_x}  \, d(x, r_Q)   \\
        &\le \sum_{ \bar{h} < \bar{l}} \binom{\bar{l}}{\bar{h}} \frac{ { (98e)^{|\bar{l} - \bar{h}|} } \, B'_{\bar{l}-\bar{h}} \, A_{m+1}^{\bar{h}} }{ \varepsilon^{|\bar{l}-\bar{h}|}} \, N_n \, (c \, \delta) = H_{\bar{l}} \, (c \, \delta) < 2^{-i-1}.
    \end{align*}
    Hence, $y''_{\bar{l}}$ is a good approximation, and the theorem is proved.
\end{proof}

\section{Open questions}
In this paper we proved the computability of Whitney Extension operators $\mathrm{WET}_m$ for $m\geq 0$. Such operators associate every closed set $F\subseteq\RR$ and every Whitney jet $(f^{(\bar{k})})_{|\bar{k}| \le m}$ of order $m$ defined of $F$ to total functions $g:\RR^n\to\RR$ such that
\begin{itemize}
    \item ${\partial_{\bar{k}}} \,g(x) = f^{(\bar{k})}(x)$ on $F$ for all $\bar{k}$ with $|\bar{k}| \le m$,
    \item $g$ is $C^\infty$ on $\RR^n \setminus F$.
\end{itemize}
This is the content of the Computable Whitney Extension Theorem \ref{Theorem_WETm}. For the case $m=0$, one obtains the First Computable Whitney Extension Theorem \ref{Theorem_FirstWET}, where only the continuity of $g$ is concerned. These results have been obtained under the assumption that closed sets are encoded via \emph{total representation}.

This leaves open at least two important questions.

The first one concerns the necessity of using total representation to encode closed sets. Although this representation is very natural, in literature closed sets are often represented  by negative information only. We conjecture that negative information does not suffice to prove computability, at least for the type of extension that we have investigated. If so, the next goal would be determining the Weihrauch degree of the corresponding operators $\mathrm{WET}_m^-$, whose definition is obtained by allowing only negative information for closed domains.

The second question concerns the possibility of proving computability for
operators WET having as input infinite Whitney jets   $(f^{(\bar{k})})_{|\bar{k}|\in\NN}$ (with closed domains encoded by total or negative representation), and as output total continuous functions $g:\RR^n\to\RR$ such that ${\partial_{\bar{k}}} \,g(x) = f^{(\bar{k})}(x)$ on $F$ for all $\bar{k}$. The extension originally defined by Hassler Whitney covered also this case.

The investigation of these open problems will be the subject for future work.

\appendix
\section{Proof of Proposition \ref{Prop_Derivativesvarphi}}\label{Appendix_A}

In order to prove Proposition \ref{Prop_Derivativesvarphi} we need a general fact about the derivatives of a quotient.

\begin{prop}\label{Prop_DerivativesQuotient}
%Let $f,g : \RR^n \to \RR$ be smooth functions. Then, for every $\bar{k}$,
%\begin{align}\label{Align_DerivativesQuotients}
%    \frac{\partial^{\bar{k}} }{\partial x^{\bar{k}}} \frac{f}{g} = \frac{Q_{\bar{k}} \left( \left\{ \frac{\partial^{\bar{l}} \, f }{\partial x^{\bar{l}}}, \frac{\partial^{\bar{l}} \, g }{\partial x^{\bar{l}}} \mid \bar{l} \le \bar{k} \right\} \right) }{g^{|\bar{k}|+1}},
%\end{align}
%where $Q_{\bar{k}}$ is a polynomial with $l_{\bar{k}}$ terms, integer coefficients $(q_j)_{j < l_{\bar{k}}}$, and each term has $\mathrm{deg}^*$ equal to $|\bar{k}|$.
%The sequence of coefficients $(q_j)_{j < l_{\bar{k}}}$ can be effectively computed from $\bar{k}$.
%
Let $u,v : \RR^n \to \RR$ be $C^m(\RR)$ functions and, for every $\bar{k}$ with $|\bar{k}| \le m$, let $\bar{l}_0, \dots, \bar{l}_{t_{\bar{k}}-1}$ list in lexicographical order all $\bar{l} \le \bar{k}$. Then the derivative ${\partial_{\bar{k}} } \frac{u}{v}$ can be expressed as
\begin{align}\label{Align_DerivativesQuotients}
\frac{\sum_{j < t_{\bar{k}}} W_j}{v^{|\bar{k}|+1}},
\end{align}
where each $W_j$ is of the form
\begin{align*}
    W_j = r_j \, ( \partial_{\bar{l}_j} u )%^{n_{j,\bar{l}}}
    \prod_{\bar{l} \le \bar{k}}  ( {\partial_{\bar{l}} \, v} )^{m_{j,\bar{l}}},
\end{align*}
with $r_j \in \ZZ$, $m_{j,\bar{l}} \in \NN$ and $|\bar{l}_j| + \sum_{\bar{l} \le \bar{k}} m_{j,\bar{l}} \, |\bar{l}| = |\bar{k}|$.

Moreover $t_{\bar{k}}$ and the sequences $(r_j)_{j < t_{\bar{k}}}$, and $(m_{j,\bar{l}})_{j < t_{\bar{k}}, \bar{l} \leq \bar{k}}$ can be effectively computed from $\bar{k}$.
\end{prop}
\begin{proof}
By induction on $|\bar{k}|$.
The case $|\bar{k}| = 0$ is obvious. 	

For the inductive step, suppose that (\ref{Align_DerivativesQuotients}) holds for $\bar{k}$.
We show that it holds also for $\bar{k} + e_i$ for every $i$ with $1 \le i \le n$. We have
\begin{align*}
     \partial_i \, \partial_{\bar{k}} \, \frac{u}{v} &= \sum_{j < t_{\bar{k}}} {\partial_i } \frac{W_j}{v^{|\bar{k}|+1}} =
      \sum_{j < t_{\bar{k}}} \frac{ {\partial_i }  W_j v^{|\bar{k}|+1} - W_j \, {\partial_i } (v^{|\bar{k}|+1}) }{v^{2(|\bar{k}|+1)}}  \\
      &= \sum_{j < t_{\bar{k}}} \frac{ {\partial_i }  W_j \, v^{|\bar{k}|+1} - (|\bar{k}|+1) W_j v^{|\bar{k}|} {\partial_i  v} }{v^{2(|\bar{k}|+1)}} = \sum_{j < t_{\bar{k}}} \frac{ {\partial_i }  W_j \, v - (|\bar{k}|+1) W_j \, \partial_i v }{v^{|\bar{k}|+2}}.
\end{align*}
Since
\begin{align*}
    {\partial_i } \, W_j & = r_j \, {\partial_i } \, \Big[ ( {\partial_{\bar{l}_j} u} ) \prod_{\bar{l} \le \bar{k}}  ( {\partial_{\bar{l}} v} )^{m_{j,\bar{l}}} \Big] = r_j \, ( {\partial_{\bar{l}_j + e_i}  u} ) \prod_{\bar{l} \le \bar{k}}  ( {\partial_{\bar{l}} \, v} )^{m_{j,\bar{l}}} \, + \\
     & \quad + r_j \, {\partial_{\bar{l}_j} u}  \sum_{\bar{h} \le \bar{k}} m_{j,\bar{h}} \, {\partial_{\bar{h} + e_i} v} \, ( {\partial_{\bar{h} } v} )^{m_{j,\bar{h}} - 1 }  \prod_{\bar{l} \le \bar{k}, \bar{l} \ne \bar{h} }  ( {\partial_{\bar{l}} v} )^{m_{j,\bar{l}} }
\end{align*}
all terms in the numerator are of the prescribed form.
\end{proof}

\begin{proof}[Proof of Proposition \ref{Prop_Derivativesvarphi}]
  Using (\ref{Align_varphi}), the definition of $\varphi_0$, Lemma \ref{Lemma_Derivativesh} and Notation \ref{Def_nInteger},
    we have
    \begin{align*}
        \left| \partial_{\bar{k}} \, \varphi_Q(x) \right| &=
        \frac{1}{l_Q^{|\bar{k}|}} \left| \partial_{\bar{k}}\, \varphi_0 \left( \frac{x - c_Q}{l_Q} \right) \right| \\
        &= \frac{1}{l_Q^{|\bar{k}|}} \left| \frac{\partial^{k_1}}{\partial x_1^{k_1}} \, \nu \left( \frac{x_1 - (c_Q)_1}{l_Q} \right) \dots \frac{\partial^{k_n}}{\partial x_n^{k_n}} \, \nu \left( \frac{x_n - (c_Q)_n}{l_Q} \right) \right| \\
        &\le \frac{1}{l_Q^{|\bar{k}|}} B_{k_1} \left(\frac{2}{\varepsilon}\right)^{k_1}  \dots  B_{k_n}\left( \frac{2}{\varepsilon} \right) ^{k_n} = B_{\bar{k}} \left(\frac{2}{ \varepsilon \, l_Q} \right) ^{|\bar{k}| }.
    \end{align*}

    For the second inequality the interesting case is when $x \in Q^*$. First, observe that using (\ref{Align_varphistar}) and Proposition \ref{Prop_DerivativesQuotient}, we have
    \begin{align*}
        \partial_{\bar{k}} \, \varphi^*_Q(x) = \frac{\sum_{j < t_{\bar{k}}} W_j}{\Phi(x)^{|\bar{k}| + 1}}.
    \end{align*}
Since each $W_j$ depends on the partial derivatives of $\Phi$ we need estimates on the absolute values of the latter. By the proof of Proposition \ref{Prop_Gx} we have that for every $Q' \in \FF_x$, $ \frac{1}{7}\, d(x,F) < \diam(Q') < 3\, d(x,F)$.
    Therefore, for every $Q' \in \FF_x$ we have
    \begin{align*}
        \diam(Q') > \frac{\diam(Q)}{ 21 }.
    \end{align*}
    Using that $\diam(Q) = \sqrt{n} \,l_Q$ and Proposition \ref{Prop_Gx} %and (\ref{Align_NormVarphilQ})
    we obtain for every $\bar{l} \leq \bar{k}$
    \begin{align}\label{Align_DerivativesPhi}
        \left| \partial_{\bar{l}} \, \Phi(x) \right| \le N_n \, B_{{\bar{l}}} \left( \frac{2}{ \varepsilon} \right)^{|\bar{l}|} \left( \frac{ 21 \, \sqrt{n} } { \diam(Q) }\right)^{|\bar{l}|}.
    \end{align}

    Now, consider a term $W_j$; by Proposition \ref{Prop_DerivativesQuotient}, we have
    \begin{align*}
        |W_j| &= |r_j| \, | {\partial_{ \bar{l}_j} }\, \varphi_Q(x) | \prod_{\bar{l} \le \bar{k}}  \left| {\partial_{\bar{l}} } \, \Phi(x) \right|^{m_{j,\bar{l}}} \le  K_j \, \left( \frac{2}{\varepsilon \diam(Q)} \right) ^{|\bar{k}|}
    \end{align*}
    where $K_j$ can be computed from $|r_j|$, $n$, $\left\{ B_{\bar{l}} \mid \bar{l} \le \bar{k} \right\}$ and $N_n$.

    Using again (\ref{Align_DefBbark}) and $\Phi(x) \ge 1$ (by Proposition \ref{Prop_FisCovering}) we obtain
    \begin{align*}
        \left| {\partial_{\bar{k}}} \, \varphi^*_Q(x) \right| \le \sum_{j < t_{\bar{k}}} { K_j \, \left( \frac2{\varepsilon \diam(Q)} \right)^{|\bar{k}|} }.
        % \left| Q_{\bar{k}}\left( \left\{ \frac{\partial^{\bar{l}}}{\partial x^{\bar{l}}} \varphi_Q(x), \frac{\partial^{\bar{l}}}{\partial x^{\bar{l}}} \Phi(x) \mid \bar{l} \le \bar{k} \right\} \right) \right| \le \\
    \end{align*}
    Setting $B'_{\bar{k}} := \sum_{j < t_{\bar{k}}} K_j$ we have the thesis.
\end{proof}

\end{document}